\documentclass[11pt]{amsart}

\usepackage{amsmath, amsthm, amssymb, amscd, amsfonts}
\usepackage{fontsmpl}
\usepackage{fontsmpl,layout}
\usepackage[all]{xy}             
\usepackage[active]{srcltx}   
\usepackage{color}

\newtheorem{lema}{Lemma}[section]
\newtheorem{theorem}[lema]{Theorem}
\newtheorem{cor}[lema]{Corollary}

\newtheorem{prop}[lema]{Proposition}

\theoremstyle{definition}
\newtheorem{definition}[lema]{Definition}

\theoremstyle{remark}
\newtheorem{obs}[lema]{Remark}

\newcommand\id{\operatorname{id}}
\newcommand\ad{\operatorname{ad}}

\newcommand\gr{\operatorname{gr}}
\newcommand\co{\operatorname{co}}

\newcommand\ord{\operatorname{ord}}

\newcommand\End{\operatorname{End}}
\newcommand\Hom{\operatorname{Hom}}

\newcommand\Aut{\operatorname{Aut}}

\newcommand\Ker{\operatorname{Ker}}

\newcommand\sg{\operatorname{sg}}

\newcommand{\eps}{\varepsilon}
\newcommand{\ot}{\otimes}
\newcommand{\com}{\Delta}

\newcommand\toba{{\mathfrak B }}

\newcommand{\trid}{\triangleright}

\newcommand{\I}{{\mathbb I}}

\newcommand{\U}{{\mathcal U}}

\newcommand{\D}{{\mathcal D}}

\newcommand{\Oc}{{\mathcal O}}
\newcommand{\ydh}{{}^H_H\mathcal{YD}}
\newcommand{\ydsigmah}{{}^{H_{\sigma}}_{H_{\sigma}}\mathcal{YD}}
\newcommand{\ydgamma}{{}^\Gamma_\Gamma\mathcal{YD}}

\newcommand{\ydsn}{{}^{\s_{n}}_{\s_{n}}\mathcal{YD}}

\def\NN{\mathbb{N}}
\def\ZZ {\mathbb{Z}}

\def\SS{\mathcal{S}}
\def\S{\mathbb{S}}
\def\DD{\mathbb{D}}
\def\s{\mathbb{S}}

\def\k{\mathbb{C}}

\def\Z{\mathcal{Z}}

\def\D{\mathcal{D}}

\def\mH{\mathcal{H}}
\def\mQ{\mathcal{Q}}

\def\lieg{\mathfrak{g}}

\newcommand{\GL}{\mathbf{GL}}

\def\QEA{U_{\textbf{q}}(\lieg_{A})}

\def\Ga{\Gamma}

\def\pf{\begin{proof}}
\def\epf{\end{proof}}

\theoremstyle{plain}

\begin{document}




\title[Multiparameter quantum groups, bosonizations and cocycles]
{Multiparameter quantum groups,\\ 
bosonizations and cocycle deformations}

\author[g. a. garc\'\i a]
{Gast\' on Andr\' es Garc\'\i a}

\address{\newline\noindent Departamento de Matem\'atica, 
Facultad de Ciencias Exactas,
Universidad Nacional de La Plata.
CONI\-CET.
Casilla de correo 172, (1900) La Plata,
Argentina.
}
\email{ggarcia@mate.unlp.edu.ar}

\thanks{\noindent 2010 \emph{Mathematics Subject Classification.}
16T05, 17B37, 17B35, 81R50. \newline \emph{Keywords:} quantum groups, quantized
enveloping algebras, quantized coordinate algebras. \newline
Partially supported by CONICET, ANPCyT, Secyt (Argentina) and GNSAGA (Italy).}

\date{\today.}

\begin{abstract} We 
describe quantum groups given by multiparametric deformations
of enveloping algebras of Kac-Moody algebras as 
a family of pointed Hopf algebras introduced by Andruskiewitsch and
Schneider associated to a generalized Cartan matrix. 
We show that under some assumptions, these 
Hopf algebras
depend only on one parameter on each connected component of the
Dynkin diagram, up to a cocycle
deformation. In particular, we obtain in this way a known result of Hu, Pei and Rosso.
\end{abstract}

\maketitle



\section{Introduction}\label{sec:intro}
Let $\Bbbk$ be an algebraically closed field of characteristic zero and  $\pi: A \to H$ a Hopf 
algebra projection that admits a Hopf algebra section $\iota: H \to A$. 
Then $A \simeq R\#H$, the 
Radford-Majid product or \textit{bosonization} of $R$ over $H$, where $R= A^{\co \pi}$ is a 
braided Hopf algebra in the category of Yetter-Drinfeld modules $\ydh$ over $H$. 
Conversely, given a braided Hopf algebra $R$ in $\ydh$, its bosonization $R\#H$ is
an ordinary Hopf algebra with a projection to $H$, see
\cite{R2}.
This fact play a crucial role in the classification of finite-dimensional
pointed Hopf algebras \cite{AS4} and in the description of quantum groups
arising from deformation of enveloping algebras of semisimple Lie algebras, in particular
their quantum Borel subalgebras.

If $A$ is a Hopf algebra such that its coradical $A_{0}$ is a Hopf subalgebra, then 
the graded object $\gr A$ associated to the coradical filtration is again a Hopf 
algebra and $\gr A \simeq R\# H$. A key point in the lifting method to classify
pointed Hopf algebras, where $A_{0}$ is a group algebra, relies on the description
of $R$ as a Nichols algebra, a graded connected braided Hopf algebra generated
by the elements of degree one. If $\mathcal{A}$ is a pointed Hopf algebra
such that $\gr \mathcal{A} \simeq R\#H$, it is said that $\mathcal{A}$ is a lifting
of $R$ over $H$. This procedure can be generalized for other types of filtrations, 
such as the standard filtration, see \cite{AC}.

An open question concerning this problem is whether all possible liftings of a
bosonization $R\#H$ can be obtained by a (left) 2-cocycle deformation 
on the multiplication. Positive answers were obtained for $H=\k  \Ga$
a group algebra by different authors
using different methods, among them,
\cite{GrM}, \cite{Mk}, \cite{Mk2} 
for the case $\Ga$ is abelian and $R$ is a quantum linear space, 
\cite{GIM} for $\Ga=\S_{n}$ a symmetric group, \cite{GM} for $\Ga=\DD_{m}$ 
a family of dihedral
groups and \cite{GIV} for pointed Hopf algebras associated to affine racks.
Moreover, in \cite{AAGMV} a systematic procedure is described to construct 
liftings as cocycle deformations via Hopf-Galois objects. 
Variations of this problem were also studied by other authors,
see for example \cite{ABM1, ABM2} and references 
therein.

Hence, a natural question is to describe the set $\Z^{2}(A,\k)$ of Hopf
2-cocycles on a bosonization $A=R\#H$ for $H$ a Hopf algebra with bijective antipode and
$R$ a braided Hopf algebra in $\ydh$. Since $H$ is a Hopf subalgebra of $A$, the restriction
of any Hopf 2-cocycle on $A$ gives a Hopf 2-cocycle on $H$. This restriction 
admits a section that gives an injective map 
$\Z^{2}(H,\k) \hookrightarrow \Z^{2}(A,\k)$. In particular, 
any Hopf 2-cocycle $\sigma$ on $H$ defines
a Hopf 2-cocycle  $\tilde{\sigma}$ on $A$.

Assume $\Bbbk = \k$ and let $\theta$ be a positive integer. 
Let $\QEA$ be the multiparameter quantum group
associated to a generalized Cartan matrix $A$ defined in \cite{HPR}, with 
$\mathbf{q} = (q_{ij})_{1\leq i,j\leq \theta}$ and $q_{ij}\in \k^{\times}, q_{ii}\neq 1 $
for all $1\leq i,j\leq \theta$.
In order to treat these quantum groups in a unified
way, we describe them explicitly in Theorem \ref{thm:qersg=ured} 
as a family of \textit{reductive} pointed Hopf algebras given in \cite{ARS} which
are associated to the generalized Cartan matrix $A$. 
A similar description
is given in \cite{H} using bicharacters and the Lusztig isomorphisms are described. 
This allows us to 
study them as quotients of bosonizations of pre-Nichols algebras. 

Using some results on cocyles on bosonizations, we show 
in Theorem \ref{thm:cocycle-def-pre-Nichols} that if $q_{ii}$ is a positive real number for all
$1\leq i \leq \theta$,
these families of pointed Hopf algebras depend only on one parameter on each connected component of the
Dynkin diagram associated to $A$, up to a cocycle
deformation.
This relation was previously described in \cite{AS3}, \cite{Ro} 
as a twist-equivalence between the 
matrices associated to the braiding; here we re-interpret them as 2-cocycle deformations. 

As a consequence, we show that the multiparameter quantum groups 
$\QEA$ are cocycle 
deformations of quantum groups
that depend only on one parameter on each connected 
component of the Dynking diagram, obtaining in case $\mathfrak{g}_{A}$ is simple, 
a result of Hu, Pei and Rosso \cite[Thm. 28]{HPR}, see
Theorem \ref{thm:qersg=ured} and Corollaries \ref{cor:multi-one-cocycle}, 
\ref{cor:multi-two-cocycle}.

The paper is organized as follows. In Section \ref{sec:prelim} we fix notation and recall 
some known facts on Hopf 2-cocycles, Yetter-Drinfeld modules, Nichols algebras and 
bosonizations. In particular, in Subsection \ref{sec:cocy-def-bos} we treat cocycle deformations
of Hopf algebras given by bosonizations $A=R\#H$ and describe the 
relation between $\Z^{2}(A,\k)$ and $\Z^{2}(H,\k)$.
At the end of this section we provide an example on finite-dimentional pointed 
Hopf algebras over symmetric groups.

In Section \ref{sec:multi-par-def} we first introduced the 
family of pointed Hopf algebras given in
\cite{ARS} and show that in case the braiding is positive, these Hopf algebras 
depend only on one parameter for each connected component
of the Dynkin diagram associated to them, up to a cocycle deformation. Then 
we explicitly describe 
the multiparameter quantum groups
given in \cite{HPR} as a family of these pointed Hopf algebras and apply to them
the results on cocycle deformations.


\section{Preliminaries}\label{sec:prelim}
In this section we fix notation and recall some definitions and known results
that are used along the paper.
 
\subsection{Conventions}
We work over the field $\k$ of complex numbers 
and by $\k^{\times}$ we denote the group of units.
If $\Ga$ is a group, we denote by
$\widehat{\Ga}$ the character group. 
By convention, $\NN = \{0,1,\ldots\}$. If $A$ is an algebra and $g \in A$
is invertible, then $g\trid  a = gag^{-1}$, $a \in A$, 
denotes the inner automorphism defined by $g$.

Our references for the theory of Hopf algebras are \cite{Mo},
\cite{Sw} and  \cite{R2}. We use standard notation for Hopf algebras; 
the comultiplication is denoted $\com$
 and the antipode $\SS$.
The left adjoint representation
of $H$ on itself is the algebra map 
$\ad : H \to \End (H)$, $\ad_{l} x(y) = x_{(1)} y\SS(x_{(2)})$, 
$x, y \in H$; we shall
write $\ad$ for $\ad_{l}$, omitting the subscript $l$ 
unless strictly needed. There is also a right adjoint action
given by $\ad_{r} x(y) = \SS(x_{(1)})yx_{(2)}$. 
Note that both $\ad_{l}$ and $\ad_{r}$ are multiplicative. 
The set of group-like elements of a
coalgebra $C$ is denoted by $G(C)$. We also denote by $C^{+} = \Ker
\eps$ the augmentation ideal of $C$, where $\eps: C\to \k$ is the
counit of $C$. Let $g,h\in G(H)$, the set of 
$(g,h)$-primitive elements is given by 
$P_{g,h}(H)=\{x\in H:\ \com(x)= x\ot g + h\ot x\}$.
We call $P_{1,1}(H)=P(H)$ the set of primitive elements.

Let $A\xrightarrow{\pi} H$ be a Hopf algebra map,
then 
$$ A^{\co H}= A^{\co \pi} =\{a \in A|\ (\id\otimes \pi)\Delta
(a) = a\otimes 1\},$$ 
denotes the subalgebra of right coinvariants and
$ ^{\co H}A= \ ^{\co \pi}A$
 denotes the subalgebra of left coinvariants.

 For $n > 0$ and $q\in \k ^{\times}$, $q\neq 1$, define
\begin{align*}
(n)_{q} &= \frac{q^{n}-1}{q-1}=q^{n-1} + \cdots +q + 1,\\
(n)_{q}! & =  (n)_{q}(n-1)_{q}\cdots (2)_{q}(1)_{q}\text{ and }(0)_{q}=1,\\
\binom{n}{k}_{q}&= \frac{(n)_{q}}{(k)_{q}(n-k)_{q}}.
\end{align*}
It is well-known that 
\begin{equation}\label{eq:qbinom}
\binom{n}{k}_{q} = q^{k}\binom{n-1}{k}_{q} + \binom{n-1}{k-1}_{q} =
\binom{n-1}{k}_{q} + q^{n-k}\binom{n-1}{k-1}_{q}.
\end{equation}

A \textit{braided vector space} is a pair $(V,c)$ with $V$ a vector space and 
$c\in \Aut(V\ot V)$ that satisfies the braid equation, that is
$(c\ot \id)(\id\ot c)(c\ot \id) = (\id\ot c)(c\ot \id)(\id\ot c) \in \End(V\ot V\ot V)$.

Let $(\mathcal{C}, \ot,a,l,r,\mathbf{1})$ be a monoidal category and denote 
by $\tau:\mathcal{C} \times \mathcal{C} \to \mathcal{C}\times \mathcal{C} $
the flip functor given by $\tau(X,Y) = (Y,X)$ for all $X,Y\in \mathcal{C}$.
A \textit{braiding} on $\mathcal{C}$ is a natural isomorphism $c: \ot \to \ot \tau$ 
that satisfies the hexagon axiom for any $U,V,W \in \mathcal{C}$:
\begin{align*}
a_{V,W,U}\, c_{U,V\ot W}\, a_{U,V,W} &= (\id_{V}\ot c_{U,W})\, a_{V,U,W}\,
(c_{U,V}\ot \id_{W}),\\
a_{W,U,V}^{-1}\,c_{U\ot V, W}\, a_{U,V,W}^{-1} &= (c_{U,W}\ot 
\id_{V})\, a_{U,W,V}^{-1}\,
(\id_{U}\ot c_{V,W}).
\end{align*}
A \textit{braided monoidal category} is a pair $(\mathcal{C},c)$ where $\mathcal{C}$
is a monoidal category and $c$ is a braiding on $\mathcal{C}$, see \cite[Ch. XIII]{kassel}
for details. If $\mathcal{C}$ is strict, the equalities above are equivalent to
$$
c_{U,V\ot W} = (\id_{V}\ot c_{U,W})\,
(c_{U,V}\ot \id_{W}),\qquad 
c_{U\ot V, W} = (c_{U,W}\ot 
\id_{V})\,
(\id_{U}\ot c_{V,W}).
$$
In particular, if $V\in \mathcal{C}$ then $(V,c_{V,V})$ is a braided vector space.

\subsection{Deforming cocycles}\label{subsubsec:def-cocycles}
Let $A$ be a Hopf algebra.
Recall that a convolution invertible linear map 
$\sigma $ in $\Hom_{\k }(A\ot A, \k )$
is a  
\textit{normalized Hopf 2-cocycle} if 
$$ \sigma(b_{(1)},c_{(1)})\sigma(a,b_{(2)}c_{(2)}) =
\sigma(a_{(1)},b_{(1)})\sigma(a_{(2)}b_{(2)},c)  $$
and $\sigma (a,1) = \eps(a) = \sigma(1,a)$ for all $a,b,c \in A$,
see \cite[Sec. 7.1]{Mo}.

Using a $2$-cocycle $\sigma$ it is possible to define 
a new algebra structure on $A$ by deforming the multiplication,
which we denote by $ A_{\sigma} $. Moreover, $A_{\sigma}$ is indeed
a Hopf algebra with 
$A = A_{\sigma}$ as coalgebras,  
deformed multiplication
$m_{\sigma} = \sigma * m * \sigma^{-1} : A \ot A \to A$
given by $$m_{\sigma}(a,b) = a\cdot_{\sigma}b = \sigma(a_{(1)},b_{(1)})a_{(2)}b_{(2)}
\sigma^{-1}(a_{(3)},b_{(3)})\qquad\text{ for all }a,b\in A,$$ 
and antipode
$\SS_{\sigma} = \sigma * \SS * \sigma^{-1} : A \to A$ given by (see \cite{doi} for details)
$$\SS_{\sigma}(a)=\sigma(a_{(1)},\SS(a_{(2)}))\SS(a_{(3)})
\sigma^{-1}(\SS(a_{(4)}),a_{(5)})\qquad\text{ for all }a\in A.$$

We denote by $\Z^{2}(A, \k )$ the set of normalized Hopf 2-cocycles on $A$.
Let $\tau,\sigma :A\ot A \to \k $ be two linear maps. We denote
by $\tau*\sigma: A\ot A \to \k $ the linear map given by the convolution, that is
$(\tau*\sigma)(a,b)= \tau(a_{(1)},b_{(1)})\sigma(a_{(2)},b_{(2)})$ for all $a,b\in A$.

\begin{obs}\label{rmk:twisted-mod-str-gr} 
Assume $A=\k  \Ga$, with $\Ga$ a group. 
Then a normalized \textit{Hopf} 2-cocycle on $A$ is equivalent to
a 2-cocycle $\varphi \in  \Z^{2}(\Ga, \k )$, that is a map 
$\varphi : \Ga \times \Ga \to \k ^{\times}$ such that
$$\varphi(g, h) \varphi(gh, t) = \varphi(h, t) \varphi(g, ht)
\text{ and }\varphi (g, e) = \varphi (e, g) = 1 \text{ for all }g, h, t \in \Ga.$$
\end{obs}

\subsection{Yetter-Drinfeld modules, Nichols algebras
and bosonization}
Let $H$ be a Hopf algebra with
bijective antipode. A Yetter-Drinfeld module over $H$
is a left $H$-module and a left $H$-comodule with
comodule structure denoted by $\delta: V \to H\ot V$, 
$v \mapsto v_{(-1)} \ot v_{(0)}$, such that
$$\delta(h\cdot v) =h_{(1)}v_{(-1)}S(h_{(3)})\ot h_{(2)}\cdot v_{(0)}\quad
\text{ for all }v \in V, h\in H.$$ 
Let $\ydh$
be the category of Yetter-Drinfeld modules over $H$ with $H$-linear and
$H$-colinear maps as morphisms.
The category $\ydh$ is monoidal and braided. Indeed, if $V,W \in\ydh$, 
then $V\ot W$ is the tensor product
over $\k $ with the diagonal action and coaction of $H$ and braiding
$c_{V ,W}: V \ot W \to W \ot V$, $v \ot w \mapsto v_{(-1)} \cdot w \ot v_{(0)}$
for all $v \in V, w \in W$. 

If $H=\k \Gamma$ is a group algebra of a group $\Gamma$, 
we denote this category simply by $\ydgamma$. In this case, 
$V \in \ydgamma$ corresponds to a $\Gamma$-graded
vector space $V =\bigoplus_{g\in\Gamma}V_{g}$ 
which is a left $\Gamma$-module such that each homogeneous component $V_{g}$,
$g \in \Gamma$, is stable under the action of $\Gamma$.
Here, the $\Gamma$-grading yields the left $\k \Gamma$-comodule
structure by $\delta: V \to\k \Gamma \ot V$, $\delta(v) = g \ot v$ 
if $v$ is homogeneous of degree $g \in \Gamma$.
For $V,W\in \ydgamma$, the braiding is given by
$c_{V,W}(v\ot w) = g\cdot w\ot v$, for all $v\in V_{g}$, $w\in W$ and
$g\in \Gamma$. 
If $\Ga$ is finite,
then $\ydgamma$ is a semisimple category. 

Let $V\in \ydgamma$ and $g\in \Ga$, $\chi \in \widehat{\Ga}$. We denote by
$$V_{g}^{\chi} = \{v\in V:\ \delta(v) = g\ot v, h\cdot v = \chi(h)v
\ \forall h\in \Ga\},$$ 
the Yetter-Drinfeld submodule given by the $g$-homogeneous
elements with diagonal action of $\Ga$ given by $\chi$. 
In case $\Ga$ is finite abelian, the
pairs $(g,\chi)$ with $g\in \Ga$ and $\chi \in \widehat{\Ga}$ parametrize 
the simple modules and for all $V\in \ydgamma$ we 
have that
$V = \bigoplus_{g\in \Ga, \chi\in \widehat{\Ga}} V_{g}^{\chi}$. 

Since $\ydh$ is a braided monoidal category, we may consider Hopf algebras in
$\ydh$. For $V \in \ydh$, 
the tensor algebra $T (V ) = \oplus_{n\geq 0} T^{n}(V )$
is an $\NN$-graded algebra and coalgebra in the
braided category $\ydh$ 
where the elements of $V = T (V )(1)$ are primitive. 

Let $(V,c)$ be a finite-dimensional braided vector space. We say that the braiding
$c:V\ot V \to V\ot V$ is \textit{diagonal} \cite[Def. 1.1]{AS3} if there exists
a basis $x_{1},\ldots ,x_{\theta}$ of $V$ and non-zero scalars $q_{ij}$ 
such that $c(x_{i}\ot x_{j}) = q_{ij} x_{j}\ot x_{i}$ for all $1\leq i,j\leq \theta$. 
The braiding is called \textit{generic} if it is diagonal and $q_{ii}$ is not a root of unity
for all $1\leq i\leq \theta$, and it is called \textit{positive} if it is generic and $q_{ii}$ is a 
positive real number for all $1\leq i \leq \theta$. We say that two finite-dimensional 
braided vector spaces
of diagonal type $(V,c)$ and $(W,d)$ with matrices $(q_{ij})$ and $(\hat{q}_{ij})$ 
are \textit{twist-equivalent} \cite[Def. 3.8]{AS2} 
if $\dim V = \dim W$, $q_{ii} = \hat{q}_{ii}$ and 
$$q_{ij}q_{ji} = \hat{q}_{ij}\hat{q}_{ji} \qquad\text{ for all }1\leq i,j\leq \theta.$$

We are particularly interested in one class of braided Hopf algebras in these
categories, which turn out to be crucial in the theory, 
the (pre-) Nichols algebras.

\begin{definition} 
Let $V \in \ydh$
and $I(V ) \subseteq T (V )$ the largest $\NN$-graded ideal and coideal such
that $I \cap V = 0$. We call $\toba(V ) = T (V )/I(V )$ the {\it Nichols algebra} of 
$V$. In particular, $\toba(V ) = \bigoplus_{n\geq 0}\toba^{n}(V )$ is an
$\NN$-graded Hopf algebra in $\ydh$.
\end{definition}

Given a 
braided vector space $(V,c)$, one may construct the 
Nichols algebra $\toba(V,c) = \toba(V)$ in a way similar to the construction above,
by taking a quotient of the tensor algebra $T(V)$ by the homogeneous 
two-sided ideal given by the kernel of a homogeneous symmetrizer: 
Let $\mathbb{B}_{n}$ be the braid
group of $n$ letters.
Since $c$ satisfies the braid
equation, it induces a representation of 
$\mathbb{B}_{n}$, $\rho_{n}:\mathbb{B}_{n} \to \GL(V^{\ot n})$
for each $n \geq 2$. 
Consider the morphisms
$$Q_{n} = \sum_{\sigma \in \s_{n}}\rho_{n}(M (\sigma)) \in \End(V^{\ot n} ),$$
where $M:\s_{n} \to \mathbb{B}_{n}$ is the Matsumoto section
corresponding to the canonical projection 
$\mathbb{B}_{n}\twoheadrightarrow \s_{n}$.
Then the Nichols algebra $\toba(V)$ 
is the quotient of the tensor algebra $T (V)$ by the two-sided ideal
$\mathcal{J}=\bigoplus_{n\geq 2 } \Ker Q_{n}$. 

A \textit{pre-Nichols algebra}
is an intermediate graded braided Hopf algebra between $T (V )$ and $\toba(V)$,
see \cite{Mk, Mk2}.

Let $R$ be a Hopf algebra in $\ydh$ with multiplication $m_{R}$. 
For $x,y \in P(R)$, we define the 
braided adjoint action of $x$ on $y$ by 
$$\ad_{c}(x)(y) = m_{R}(x\ot y) - m_{R}\circ c_{R\ot R} (x\ot y)=
xy - (x_{(-1)}\cdot y) x_{(0)}.$$
This element is also called the {\it braided commutator} of $x$ and $y$. 

\subsubsection{Bosonization and Hopf algebras with a projection}
Let $R$ be a Hopf algebra in $\ydh$. The procedure to obtain a 
usual Hopf algebra from the braided Hopf algebra $R$ and $H$ is called 
bosonization or Radford-Majid product, and it is usually
denoted by $R \#H$. As a vector space, $R \# H = R\otimes H$ and the
multiplication and comultiplication are given by the smash-product and 
smash-coproduct,
respectively. That is, for all $r, s \in R$ and $g,h \in H$, we have
\begin{align*}
(r \# g)(s \#h) & = r(g_{(1)}\cdot s)\# g_{(2)}h,\\ 
\com(r \# g) & =r^{(1)} \# (r^{(2)})_{(-1)}g_{(1)} \ot 
(r^{(2)})_{(0)}\# g_{(2)},\\
\SS(r\# g) & = (1\# \SS_{H}(r_{(-1)}g))(\SS_{R}(r_{(0)})\# 1),
\end{align*}
where $\com_{R}(r) = r^{(1)}\ot r^{(2)}$ denotes the comultiplication in $R\in \ydh$ 
and $\SS_{R}$ the antipode. 
Clearly, the map $\iota: H \to R\#H$ given by $\iota(h) = 1\#h$ for all
$h\in H$ is an injective Hopf algebra map, and the map $\pi: R\#H \to H$ 
given by $\pi(r\#h) = \eps_{R}(r)h$ for all $r\in R$, $h\in H$
is a surjective Hopf algebra such that $\pi \circ \iota = \id_{H} $. 
Moreover, it holds that $R = (R\#H)^{\co \pi}$. 

Conversely, let $A$ be a Hopf algebra with bijective antipode and
$\pi: A\to H$ a Hopf algebra epimorphism admiting 
a Hopf algebra section $\iota: H\to A$ such that $\pi\circ\iota =\id_{H}$.
Then $R=A^{\co\pi}$ is a braided Hopf algebra in $\ydh$ called the 
\textit{diagram} of $A$ and $A\simeq R\# H$
as Hopf algebras. See \cite[11.6]{R2} for more details.

\subsection{On cocycle deformations and bosonizations}
\label{sec:cocy-def-bos}
In this subsection we collect some results on the construction of 
2-cocycles on bosonizations of Hopf algebras. 

Let $H$ be a Hopf algebra with bijective antipode, $R$ a braided
Hopf algebra in $\ydh$ and $A=R\#H$. 
To avoid confusion, in this section we denote  
by $\rightharpoonup: H\ot R \to R $ the action
of $H$ on $R$.

Let $\sigma\in \Z^{2}(H, \k )$. Then the map $\tilde{\sigma}: A\ot A \to \k $ given by 
$$ \tilde{\sigma}(r\#h, s\# k)= \sigma(h,k)\eps_{R}(r)\eps_{R}(s) 
\qquad \text{for all }r,s\in R,\ h,k\in H,$$
is a normalized Hopf 2-cocycle such that $\tilde{\sigma}|_{H\ot H}= \sigma$. 
Moreover, $H_{\sigma}$ is a Hopf subalgebra of $A_{\tilde{\sigma}}$ and 
the map
$\Z^{2}(H, \k ) \to \Z^{2}(A, \k )$ given by 
$\sigma \mapsto \tilde{\sigma}$ gives a section 
of the map $\Z^{2}(A, \k ) \to \Z^{2}(H, \k )$ 
induced by the restriction; in particular, it is injective. 
See e.g. \cite[Section 5]{Mk2}, \cite[Prop. 4.2]{CP}.

\begin{prop}\cite[Prop. 5.2]{Mk2}\label{prop:cocycle-lifting-bos}
Let $\sigma$ and $\tilde{\sigma}$ be as above.
Then $A_{\tilde{\sigma}} = R_{\sigma}\# H_{\sigma}$, 
where $R_{\sigma}=R$ as coalgebras, 
and the product is given by
\begin{equation}\label{eq:cocycle-lifting-bos}
a\cdot_{\sigma} b = \sigma(a_{(-1)},b_{(-1)}) a_{(0)}b_{(0)} 
\text{ for all }a,b \in R.
\end{equation}
Moreover,
$R_{\sigma} \in \ydsigmah$
with the action of $H_{\sigma}$ given by 
\begin{equation}\label{eq:act-sig-R}
h \rightharpoonup_{\sigma} a = \sigma (h_{(1)} , a_{(-1)} )(h_{(2)} 
\rightharpoonup a_{(0)} )_{(0)} 
\sigma^{-1} ((h_{(2)} \rightharpoonup a_{(0)} )_{(-1)} ,h_{(3)} )
\text{ for all }h \in H_{\sigma}, a \in R_{\sigma}.
\end{equation}\qed
\end{prop}

\begin{obs}
In case $R$ is a (pre-) Nichols algebra, by \cite{AAGMV}, \cite[2.7, 3.4]{MO}
$R_{\sigma}$ is also a
(pre-) Nichols algebra and the action is described by \eqref{eq:act-sig-R}. 
\end{obs}

\begin{obs}
Assume $H=\k  \Ga$ and $\sigma \in \Z^{2}(\Ga, \k )$. 
Let $h\in \Ga$ and $a\in R$ be a 
homogeneous element of degree $g\in \Ga$; in particular, 
$\delta(a) = g\otimes a$ and $\com_{A}(a)=a\ot 1+ g\ot a$.
Then \eqref{eq:act-sig-R} yields
\begin{equation}\label{eq:action-cocycle-skew-prim}
h\rightharpoonup_{\sigma} a =
h\cdot_{\sigma}a\cdot_{\sigma}h^{-1}= \sigma(h,g)
\sigma^{-1}(hgh^{-1},h) h\rightharpoonup a.
 \end{equation}
\end{obs}

\begin{obs}\label{rmk:def-ast}
 In \cite{AST} the authors introduced another type of cocycle
deformation on a Hopf algebra, 
which is closely related to the one given above. We describe
it shorlty.
Let $\Ga$ be an abelian group and $ A $ a Hopf algebra  
   that is  $ \Ga\times \Ga $-graded. 
   Given any $\varphi \in \Z^{2}(\Ga,\k)$,  
   define a new product on $A$ by
\begin{equation}\label{eq:prod-star-ast}
   h \mathop{*}\limits_{\hskip-1,3pt \varphi} k :=   
  {\varphi(\eta,\kappa)} \varphi\big(\eta',\kappa'\big)^{-1} \, h \cdot k
 \end{equation}
for all homogeneous  $ \, h , k \in H$  with degrees  $\big(\eta,\eta'\big), 
\big(\kappa,\kappa'\big) \in \Ga \times \Ga$. With this
multiplication, the new algebra $A^{(\varphi)}$ is a Hopf algebra 
with the same coalgebra structure and unit as  $ A $. Assume
$A = R\#\k\Ga$ is given by a bosonization over an abelian group $\Ga$. 
Then, the coaction of $\k \Ga$ on the elements of $R$ induces a $\Ga\times \Ga$
grading on $A$ with $\deg g = (g,g)$ for all $g\in \Ga$ and 
$\deg(x) = (g,1)$ if $\delta(x) = g\ot x$ with $x\in R$ a homogeneous element
and we have that $A^{(\varphi)} = A_{\tilde{\varphi}}$, where
$\tilde{\varphi}$ is the Hopf 2-cocycle on $A$ induced by $\varphi$. 
In particular, for $x,y$ homogeneous elements of $R$ of degree $g$ and $h$
respectively, we have that 
$$x\mathop{*}\limits_{\hskip-1,3pt \varphi} y =
\varphi(g,h)\varphi(1,1)^{-1}xy = \varphi(x_{(-1)},y_{(-1)}) x_{(0)}y_{(0)},$$
 which coincides with formula \eqref{eq:cocycle-lifting-bos}.
\end{obs}

\begin{obs}\label{rmk:cocyclegroupoid}
Let $A$ be a Hopf algebra and 
$\sigma \in \Z^{2}(A,\k )$, $\tau \in \Z^{2}(A_{{\sigma}}, \k )$. Then
$\tau*{\sigma} \in \Z^{2}(A,\k )$.
As is it known, the set $\Z^{2}(A,\k )$ is not a group in general.
It is the case if the cocycles are \textit{lazy}, see for example \cite{CP}.
Nevertheless, there is a groupoid 
structure on the set of all Hopf 2-cocycles.

Let $\Z$ be the groupoid whose objects are Hopf algebras $A$ and arrows labelled by
the set of 2-cocycles
$\{\alpha_{\sigma}:A \to A_{\sigma}:\ \sigma \in \Z^{2}(A,\k )\}$. The source and
target maps are given by 
$s(\alpha_{\sigma})= A, t(\alpha_{\sigma})=A_{\sigma}$ and the composition 
by $\alpha_{\tau}\circ \alpha_{\sigma} = \alpha_{\tau*\sigma}$ for 
$\sigma \in \Z^{2}(A,\k )$ and $\tau \in \Z^{2}(A_{\sigma},\k )$.

Clearly, the identity arrow is given by $\id_{A} = \alpha_{\eps_{A}}$, and 
since $(A_{\sigma})_{\sigma^{-1}} = A_{\sigma*\sigma^{-1}} = 
A = A_{\sigma^{-1}*\sigma}=
(A_{\sigma^{-1}})_{\sigma}$, each 
arrow is invertible with inverse $\alpha_{\sigma^{-1}}:A_{\sigma} \to A$.
\end{obs}

\subsubsection{An example on pointed Hopf algebras over $\s_{n}$}
 We describe now an example where two non-isomorphic families
 of finite-dimensional pointed Hopf algebras over $\s_{n}$ 
 are cocycle deformation of each other. The $2$-cocycle is
 given by the product of $2$-cocycle associated to
a Hochschild 2-cocycle on a Nichols algebra and a group 2-cocycle
 founded by Vendramin \cite{V} which serves as twisting of the 2-cocycles asociated
 to the rack of traspositions in $\s_{n}$. 
For this purpose we need to introduce first some 
terminology, see \cite[Def. 1.1]{AGr} for 
more details.

\subsubsection*{Racks and Nichols algebras}\label{subsubsec:racks}
A \emph{rack} is a pair $(X,\rhd)$,
where $X$ is a
non-empty set and $\rhd:X\times X\to X$ is a function, such that
$\phi_i=i\rhd (\cdot):X\to X$ is a bijection for all $i\in X$
satisfying that
$
i\rhd(j\rhd k)=(i\rhd j)\rhd (i\rhd k) \text{ for all }i,j,k\in X
$. A group $G$ is a rack with $x \trid y = xyx^{-1}$ for all 
$x, y \in G$. If $G=\s_{n}$, then we denote by $\Oc_j^n$
the conjugacy class of all $j$-cycles in
$\s_n$. 
 
 Let $(X,\rhd)$ be a rack. A
\textit{rack 2-cocycle} $q:X\times X\to \k ^{\times}$,
$(i,j)\mapsto q_{ij}$ is a function such that 
$$ q_{i,j\rhd k}\, q_{j,k}=q_{i\rhd j,i\rhd
k}\, q_{i,k}, \text{ for all }\,i,j,k\in X.$$
It determines a braiding $c^q$ on the
vector space $\k  X$ with basis $\{x_i\}_{i\in X}$ by $c^q(x_i\otimes
x_j)=q_{ij}x_{i\rhd j}\otimes x_i$ for all $i,j\in X$. We denote 
this braided vector space $(\k X, c^{q})$ by $M(X,q)$ and the Nichols algebra 
associated with it 
by $\toba(X,q)$.

Let $X$ be a subrack of a conjugacy class $\Oc$ in $\Ga$, $q$ a rack 
2-cocycle on $X$ 
and $\varphi \in \Z^{2} (\Ga,\k )$.
Then the map $q^{\varphi} : X\times X \to \k ^{\times}$ given by
\begin{equation}\label{eq:rack-cocycle-gr-cocycle}
q^{\varphi}_{xy} = \varphi(x, y)\varphi^{-1} (x \trid y, x) q_{xy}, \text{ for all }
x, y \in X,
\end{equation}
is a rack 2-cocycle.

If $X$ is any rack, $q$ a rack 2-cocycle on $X$ and $ \varphi: 
X \times X \to \k^{\times}$, then
define $q^{\varphi}$ by \eqref{eq:rack-cocycle-gr-cocycle}. 
It can be shown that $q^{\varphi}$ is a rack 2-cocycle if and only if
$$\varphi(x, z)\varphi(x\trid y, x \trid z)\varphi(x \trid(y \trid z), x)\varphi(y\trid z, y)
= \varphi(y, z)\varphi(x, y \trid z)\varphi(x\trid (y\trid z), x\trid y)\varphi(x\trid z, x)$$
for any $x, y, z \in X$. If $X$ is a subrack of a group $\Ga$ and $\varphi \in 
Z^{2}(\Ga, \k )$, then $\varphi|_{X\times X}$ satisfies the equation above.

\begin{definition}
Let $q,q': X\times X \to \k ^{\times}$ be rack 2-cocycles on $X$. 
We say that $q$ and $q'$ are \textit{twist equivalent} 
if there
exists $\varphi: X \times X \to \k ^{\times}$ such that 
$q' = q^{\varphi}$ as in 
\eqref{eq:rack-cocycle-gr-cocycle}. 
\end{definition}

\subsubsection*{On Nichols algebras over $\s_{n}$}
Let $X=\Oc_2^n$ be the rack of traspositions with $ n\geq 3 $ 
and consider the cocycles:
\begin{align*}
 -1&:\Oc_2^n\times \Oc_2^n\to \k ^{\times}, 
& & (j,i)\mapsto \sg(j)=-1;\\
 \chi &: \Oc_2^n\times \Oc_2^n\to 
\k ^{\times}, & & (j,i)\mapsto \chi_i(j) =
\begin{cases}
  1,  & \mbox{if} \  i=(a,b) \text{ and } j(a)<j(b), \\
  -1, & \mbox{if} \ i=(a,b) \text{ and } j(a)>j(b).
\end{cases}
\end{align*}
for all $i,j\in  \Oc_2^n$. By \cite[Ex. 6.4]{MS}, 
\cite[Thm. 6.12]{AGr2},
the Nichols algebras are
given by 
\begin{enumerate}
 \item[$(a)$] ${\toba}(\Oc_2^n,-1)$; generated
by the elements $\{x_{(\ell m)}\}_{ 1\le \ell < m \le n}$ satisfying
for all $1\le a < b < c \le n, 1\le e < f \le n, \{a,b\}\cap\{e,f\}=\emptyset$
the identities
\begin{align*}
0= x_{(ab)}^2  = x_{(ab)}x_{(ef)}+x_{(ef)} x_{(ab)} = 
x_{(ab)}x_{(bc)}+x_{(bc)} x_{(ac)}+x_{(ac)} x_{(ab)}.
\end{align*}
 \item[$(b)$] ${\toba}(\Oc_2^n,\chi)$; 
generated by the elements $\{x_{(\ell m)}\}_{ 1\le \ell < m \le n}$ satisfying
for all $1\le a < b < c \le n, 1\le e < f \le n, \{a,b\}\cap\{e,f\}=\emptyset$ the identities
\begin{align*}
0 &= x_{(ab)}^2 =  x_{(ab)}x_{(ef)} - x_{(ef)} x_{(ab)}
= x_{(ab)}x_{(bc)} - x_{(bc)} x_{(ac)} - x_{(ac)} x_{(ab)},\\
0&= 
x_{(bc)}x_{(ab)} -  x_{(ac)}x_{(bc)} - x_{(ab)}x_{(ac)}.
\end{align*}
\end{enumerate}
 
For $3\leq n \leq 5$ these Nichols algebras 
are finite-dimensional. If $n>5$ it is not known
if this is the case. It turns out that the cocycles associated to them are twist equivalent.

\begin{theorem}\cite[Thm. 3.8]{V}\label{thm:vendra-twist}
Let $n \geq 4$. The rack 2-cocycles $\chi$ and $-1$ associated to $\Oc_2^n$
are twist equivalent. \qed
\end{theorem}

\begin{obs}\label{rmk:twist-gr-cocycle}
The twist given by Theorem \ref{thm:vendra-twist} is defined using a group 2-cocycle
$\varphi \in \Z^{2}(\s_{n},\k )$. In particular,
$ -1 = \varphi(x,y)\varphi^{-1}(x\trid y, x)\chi(x,y) $
 for all $x,y \in \Oc_2^n$.
\end{obs}
 
\subsubsection*{Cocycles on pointed Hopf algebras over $\s_{n}$} 
Assume $n\geq 4$. 
Let $\varphi \in \Z^{2}(\s_{n}, \k )$ be the group 2-cocycle given in 
Remark \ref{rmk:twist-gr-cocycle}.
Denote again by $\varphi$ the associated Hopf 2-cocycle in
 $\Z^{2}(\k \s_{n},\k )$ and
by $\tilde{\varphi}\in \Z^{2}(A,\k )$ the Hopf 2-cocycle on the bosonization
$A=\toba (\Oc_2^n , \chi) \#\k  \s_{n} $. Then by 
Proposition \ref{prop:cocycle-lifting-bos}
we have that 
$$\toba ( \Oc_2^n, -1)\# \k  \s_{n} \simeq \toba (\Oc_2^n , \chi)_{\varphi}
\#\k  \s_{n} \simeq (\toba (\Oc_2^n , \chi)
\#\k  \s_{n})_{\tilde{\varphi}}. $$

Let $\Lambda, \Gamma \in \k$ and $t=(\Lambda, \Gamma) $.
Denote by $\mH(\mQ_n^{-1}[t])$ the
algebra generated by $\{a_i, h_r :
i\in\Oc_2^n,r\in\s_n\}$ satisfying the following relations for 
$r,s,j\in\s_n$ and $i\in \Oc_2^n$:
\begin{align*}
h_e=1, \quad h_rh_s=h_{rs}, \qquad 
h_j a_i=-a_{j\trid i}h_j, 
\quad a_{(12)}^2=0,\\
a_{(12)} a_{(34)}+a_{(34)} a_{(12)}=
\Lambda(1-h_{(12)}h_{(34)}),\\
 a_{(12)} a_{(23)}+a_{(23)} a_{(13)}+a_{(13)}
a_{(12)}=\Gamma(1-h_{(12)}h_{(23)}).
\end{align*}
This algebra is indeed a Hopf algebra 
with the structure determined by $h_{\sigma}$ being a grouplike
element and $a_{\sigma}$ being a $(1,h_{\sigma})$-primitive
for all $\sigma \in \Oc^{n}_{2}$. Consequently, it is  
a pointed Hopf algebra with diagram $\toba ( \Oc_2^n, -1)$, see \cite{GG}
for details. If $t= (2\lambda,3\lambda)$ with $\lambda \in \k^{\times}$, we know
that $\mH(\mQ_n^{-1}[t])$
is a cocycle deformation of the bosonization $\toba ( \Oc_2^n, -1)\# \k  \s_{n}$. 
The explicit cocycle is given in the theorem below; it was also shown in \cite{GIM}
by other methods. 

Shortly, let $X$ be a rack, $q$ a rack 2-cocycle and $\{x_{\tau}\}_{\tau \in X}$
be homogeneous elements in $V= M(X,q) \in \ydsn$. Then the linear combination
of tensor products of 
linear functionals $\delta_{\tau}$ given
by $\delta_{\tau}(x_{\mu}) = \delta_{\tau,\mu}$ for all $\mu,\tau \in X$ give rise to 
a Hochschild 2-cocycle $\eta = \sum_{\tau,\mu \in X} a_{\tau,\mu} d_{\tau}\ot \mu$ 
by defining it via
$$\eta(\toba^{m}(V)\otimes \toba^{n}(V))=0
\text{ if }(m,n)\not= (1,1).$$

If this cocycle is invariant under the action of $\s_{n}$, \textit{i.e.}
$\eta^{h}(x,y)=\eta(h_{(1)}\rightharpoonup x, h_{(2)}\rightharpoonup y) = \eta(x,y)$ for all
$x,y\in \toba(V)$ and $h\in \s_{n}$
then one may define a Hochschild 2-cocycle 
on $A= \toba(V)\# \k  \s_{n}$
by 
$$ \tilde{\eta}(x\#h, y\#k) = \eta(x, h\rightharpoonup y)\eps(k) 
\text{ for all }x,y\in \toba(V), h,k\in \s_{n}.
$$
Moreover, 
$
\sigma=e^{\tilde{\eta}}=\sum_{i=0}^\infty \frac{ \tilde{\eta}^{* i}}{i!}
\colon A\otimes A\to \k
$
is a well-defined convolution invertible map with convolution
inverse $e^{- \tilde{\eta}}$. By 
\cite[Cor. 2.6]{GM}, this map $\sigma = e^{\tilde{\eta}}$ is a Hopf 
2-cocycle. See \cite[2.1]{GM} for more details. 

\begin{theorem}\label{thm:q-1-cocycle}\cite[Thm. 4.10 $(i)$]{GM}
Let $ A = \toba( \Oc^{n}_{2},-1)\# \k  \s_{n}$ 
and $\sigma_{\lambda} = e^{\tilde{\eta}_{\lambda}}$ a Hopf $2$-cocycle 
with $\eta_{\lambda} = \frac{\lambda}{3} 
\sum_{\mu,\tau \in \Oc_2^{n}}
d_{\tau}\ot d_{\mu} $ and $\lambda \in \k $. 
Then
$ A_{\sigma_{\lambda}} \simeq  
\mH(\mQ_n^{-1}[(2\lambda,3\lambda)])$  for $n\geq 4$.\qed 
\end{theorem}

We end this section with the following result.

\begin{cor}\label{cor:fomin-kirillov-cocycle}
Let  $\sigma_{\lambda} = e^{\tilde{\eta}_{\lambda}} 
\in \Z^{2}(\toba(\Oc_2^{n},-1)\# \k  \s_{n},\k  )$ and
$\tilde{\varphi}\in  \Z^{2}(\toba(\Oc_2^n , \chi)\# 
\k  \s_{n},\k  )$ be the Hopf 2-cocycles
defined above. Then
$$\mH(\mQ_n^{-1}[(2\lambda,3\lambda)]) \simeq (\toba (\Oc_2^n , \chi) \#\k  \s_{n})
_{\sigma_{\lambda}*\tilde{\varphi}}. $$
\end{cor}

\pf
Since  $\sigma_{\lambda} = e^{\tilde{\eta}_{\lambda}}$ is a Hopf 
2-cocycle on $\toba(\Oc_2^n,-1)\# \k  \s_{n}$
and this algebra is isomorphic to 
$ (\toba (\Oc_2^n , \chi) \#\k  \s_{n})_{\tilde{\varphi}}$ for 
$\varphi \in \Z^{2}(\s_{n},\k )$, the claim follows by Remark \ref{rmk:cocyclegroupoid}.
\epf

\section{Multiparameter quantum groups, pre-Nichols algebras and 
cocycle deformations}
\label{sec:multi-par-def}
In this section we show explicitly that certain classes of
multiparameter quantum groups can be described using the theory of pointed
Hopf algebras developed by 
Andruskiewitsch and Schneider \cite{AS4}, \cite{ARS}. 

First we introduce these families of pointed Hopf algebras and show
that under some assumptions they are cocycle deformations of 
certain (one-parameter) families of
pointed Hopf algebras. 

\subsection{On pointed Hopf algebras associated to generalized Cartan matrices} 
Let $\theta$ be a positive integer and $(a_{ij})_{1\leq i,j \leq \theta}$ a generalized Cartan matrix, 
that is, $(a_{ij})_{1\leq i,j \leq \theta}$ is a matrix with integer
entries such that $a_{ii} = 2$ for all $1 \leq i \leq \theta$, 
and for all $1 \leq i, j \leq \theta$, $i \neq j$, $a_{ij} \leq 0$, 
and if $a_{ij} = 0$, then $a_{ji} = 0$.

\begin{definition}\cite[Def. 3.2]{ARS} 
A {\it reduced YD-datum of Cartan type} 
$$\D_{red} = \D(\Gamma, (L_{i})_{1\leq i\leq\theta},(K_{i})_{1\leq i\leq\theta},  
(\chi_{i})_{1\leq i\leq \theta}, (a_{ij})_{1\leq i,j \leq \theta})$$  
consists of an abelian group $\Gamma$, $K_{i}, L_{i}\in \Gamma$
and characters $\chi_{i} \in \widehat{\Gamma} 
= \Hom(\Gamma, \k^{\times})$ satisfying
for all $1\leq i,j\leq \theta$ that
\begin{align*}
K_{i}L_{i}& \neq 1 ,\qquad
q_{ij}=\chi_{j}(K_{i}) = \chi_{i}(L_{j}),\\ 
q_{ij}q_{ji}& =q_{ii}^{a_{ij}},\quad q_{ii} \neq 1,\ 0 \leq -a_{ij} < \ord(q_{ii}) \leq \infty.
\end{align*} 
A reduced YD-datum $D_{red}$ is called {\it generic} 
if $\chi_{i}(K_{i})=q_{ii}$ is not a root of unity, for all $1\leq i \leq \theta$. 
A {\it linking parameter} $\ell$ for $\D_{red}$ is a family 
$\ell = (\ell_{i})_{1\leq i\leq \theta}$ of non-zero elements in $\k$.
\end{definition}

Let $\I = \{1, 2, \ldots , \theta\}$. We have an equivalence
relation on $\I$: for $i\neq j \in \I$ we say that
$i\sim j$ if and only if there are $i_{1},\ldots, i_{t} \in \I$ with $t\geq 2$,
$i_{1}=i, i_{t}=j$ and $q_{i_{k},i_{k+1}}q_{i_{k+1},i_{k}}\neq 1$ for all
$1\leq k<t$. We denote by $\mathcal{X}$ the set of equivalence classes.
This equivalence can be described as usual in terms of the Cartan matrix. 
Indeed, for
all $1 \leq i, j \leq \theta$, $i \sim j$ if and only if there are 
$i_{1},\ldots , i_{t} \in \I$, $t \geq 2$ with $i_{1} = i, i_{t} = j$, and
$a_{i_{k},i_{k+1}}\neq 0$
for all $1 \leq k <t$.

A reduced YD-datum of Cartan type is said of \textit{DJ-type} (Drinfeld-Jimbo type)
if the Cartan matrix is symmetrizable, \textit{i.e.} there exist $d_{i}$ relatively prime positive integers such that
$d_{i}a_{ij} = d_{j}a_{ji}$ for all $1\leq i,j\leq \theta$, and for all $I\in \mathcal{X}$ there exists
$q_{I}\in \Bbbk$
such that 
$$ q_{ij} = q_{I}^{d_{i}a_{ij}} \text{ for all }i\in I, 1\leq j\leq \theta.$$
In particular, $q_{ij}=1$ if $a_{ij} = 0$. We denote this datum by
$$\D_{q} = (\Gamma, (L_{i})_{1\leq i\leq\theta},(K_{i})_{1\leq i\leq\theta},  
(q_{I})_{I\in \mathcal{X}}, (a_{ij})_{1\leq i,j \leq \theta}).$$

\begin{definition}\label{def:pre-nichols-datum}\cite[Def. 2.4]{AS4}
Let $\D_{red}= (\Gamma, (L_{i})_{1\leq i\leq\theta} , (K_{i})_{1\leq i\leq\theta},
(\chi_{i})_{1\leq i\leq \theta}, (a_{ij})_{1\leq i,j \leq \theta}) $ be a reduced YD-datum of Cartan type and 
$\ell = (\ell_{i})_{1\leq i\leq \theta}$ a linking parameter.
$\widetilde{\U}(\D_{red},\ell)= \widetilde{\U}$ is the algebra generated by the elements 
$g \in \Gamma$, $x_{i},y_{i}$ with $1\leq i \leq \theta$ satisfying the relations
\begin{align*}
\quad & g^{\pm 1}h^{\pm 1}  
= h^{\pm 1}g^{\pm 1}, 
\quad g^{\pm 1}g^{\mp 1} 
= 1,\\
\quad & g x_{i} g^{-1}= \chi_{i}(g) x_{i},\qquad
g y_{i} g^{-1}= \chi_{i}(g)^{-1} y_{i},\\
\qquad & 
x_{i}y_{j} - \chi^{-1}_{j}(K_{i})y_{j} x_{i} = 
- \delta_{ij} \ell_{i} (K_{i}L_{i} - 1),  \\
\qquad &\ad_{c} (x_{i})^{1-a_{ij}} (x_{j})=0,
\qquad \ad_{c} (y_{i})^{1-a_{ij}} (y_{j}) =0, 
\end{align*}
for all $g,h\in \Gamma$, $1\leq i,j\leq\theta$, where
\begin{align*}
\ad_{c} (x_{i})(x_{j}) &=x_{i} x_{j} -  \chi_{j}(K_{i}) x_{j}x_{i} =  
x_{i}x_{j} - q_{ij}x_{j}x_{i},\ 1 \leq i, j \leq \theta,\\
\ad_{c} (y_{i})(y_{j}) &=y_{i}y_{j}  
-\chi^{-1}_{j}(L_{i})y_{j}y_{i} =
y_{i}y_{j} - q_{ji}^{-1}y_{j} y_{i},\ 1 \leq i, j \leq \theta.
\end{align*} 
\end{definition}
The algebra $\widetilde{\U}$ is a Hopf algebra with its structure determined by
$g \in \Gamma$ being grouplike, $x_{i}$ being $(1,K_{i})$-primitive and 
$y_{i}$ being $(1,L_{i})$-primitive for all $1\leq i \leq \theta$. In particular, it is a pointed Hopf algebra
with $G(\widetilde{\U}) = \Ga$.

\begin{obs}
The data $(\Gamma, (L_{i})_{1\leq i\leq\theta} , (K_{i})_{1\leq i\leq\theta},
(\chi_{i})_{1\leq i\leq \theta}) $ is called simply a YD-reduced datum.  Let
$V, W$ be vector spaces with basis 
$\{x_{i}\}_{1 \leq i  \leq {\theta}}$ and 
$\{y_{i}\}_{1 \leq i  \leq {\theta}}$, respectively.
Then any YD-reduced datum defines on $V$ and $W$ a Yetter-Drinfeld module structure over 
$\k \Ga$ given by $x_{i} \in V_{K_{i}}^{\chi_{i}}$ and 
$y_{i} \in W_{L_{i}}^{\chi^{-1}_{i}}$ for all $1\leq i\leq \theta$, that is,
\begin{align*}
\delta(x_{i}) & = K_{i} \ot x_{i}, & g\cdot x_{i} & = \chi_{i}(g) x_{i},\\
\delta(y_{i}) & = L_{i} \ot y_{i}, & g\cdot y_{i} & = \chi_{i}^{-1}(g) y_{i},
\end{align*}
for all $g\in \Ga$. 
The braidings $c_{V} = c_{V,V}$ of $V$ and 
$c_{W}$ of $W$ are given by
\begin{align*}
c_{V}(x_{i} \ot x_{j}) &=K_{i}\cdot x_{j}\ot x_{i} 
= \chi_{j}(K_{i}) x_{j} \ot x_{i} =  q_{ij}x_{j} \ot x_{i},\ 1 \leq i, j \leq \theta,\\
c_{W}(y_{i} \ot y_{j}) &=L_{i}\cdot y_{j}\ot y_{i}  
= \chi^{-1}_{j}(L_{i}) y_{j} \ot y_{i} =
q_{ji}^{-1}y_{j} \ot y_{i},\ 1 \leq i, j \leq \theta,
\end{align*} 
in particular, they are of diagonal type, and the corresponding adjoint actions are given by
\begin{align*}
\ad_{c} (x_{i})(x_{j}) &=x_{i} x_{j} -  \chi_{j}(K_{i}) x_{j}x_{i} =  
x_{i}x_{j} - q_{ij}x_{j}x_{i},\ 1 \leq i, j \leq \theta,\\
\ad_{c} (y_{i})(y_{j}) &= y_{i}y_{j}  
-\chi^{-1}_{j}(L_{i}) y_{j}y_{i} =
y_{i}y_{j} - q_{ji}^{-1}y_{j} y_{i},\ 1 \leq i, j \leq \theta.
\end{align*} 

The pre-Nichols algebras $R(\D)$, $R(\D, V)$ and $R(\D, W)$ associated to the reduced YD-datum
described above are given by the quotient (braided) Hopf algebras
\begin{align*}
R(\D)& = T(V\oplus W) / (\ad_{c} (x_{i})^{1-a_{ij}} (x_{j}),\ 
\ad_{c} (y_{i})^{1-a_{ij}} (y_{j}), 1\leq i\neq j \leq \theta),\\
R(\D, V)& = T(V) / (\ad_{c} (x_{i})^{1-a_{ij}} (x_{j}), 1\leq i\neq j \leq \theta),\\
R(\D, W)& = T(W) / (\ad_{c} (y_{i})^{1-a_{ij}} (y_{j}), 1\leq i\neq j \leq \theta).
\end{align*}
Since $c_{W,V} c_{V,W} = \id $, we have that $R(\D) \simeq R(\D,V)\otimes R(\D,W)$,
see \cite{Mk2}.
By abuse of notation, we denote the
images of the elements $x_{i}, y_{j}$ in $R(\D)$ again by $x_{i}, y_{j}$.
It is well-known that the elements $\ad_{c} (x_{i})^{1-a_{ij}} (x_{j}) $,
$1\leq i\neq j \leq \theta $
are primitive in the free algebra $T(V)$ 
(see for example \cite[A.1]{AS1}),
hence they generate a Hopf ideal. 

It follows that $\widetilde{\U}(\D_{red}, \ell)$ is the Hopf algebra given by
the quotient of the bosonization $R(\D)\# \k  \Ga$ 
modulo the
ideal generated by the elements
$$x_{i}y_{j} - \chi^{-1}_{j}(K_{i})y_{j} x_{i} 
- \delta_{ij} \ell_{i} (K_{i}L_{i} - 1)\qquad\text{ for all } 1\leq i,j\leq\theta,$$
where we identify $x_{i} = x_{i}\# 1$, $y_{i}=y_{i}\#1$ and 
$K_{i}= 1\# K_{i}$, $L_{i}= 1\# L_{i}$ for all $1\leq i \leq \theta$.
\end{obs}

\begin{obs}
In case the braiding is positive and generic, the pre-Nichols algebras $R(\D,V)$
and $R(\D,W)$
coincide with the Nichols algebras $\toba(V)$ and $\toba(W)$ respectively, 
that is, the ideals $I(V)$ and $I(W)$
are generated by the {\it quantum Serre relations} 
$\ad_{c} (x_{i})^{1-a_{ij}} (x_{j}) $,
$\ad_{c} (y_{i})^{1-a_{ij}} (y_{j})$ associated to the 
braided commutators, see \cite[Thm. 4.3]{AS3} and references therein. These Serre
relations are not enough to define the ideal $I(V\oplus W)$, see Remark \ref{rmk:ideal-nich-v+w}.
\end{obs}

\begin{obs}\label{rmk:ideal-nich-v+w} 
In case the braiding is positive and generic,
the ideal 
$I(V\oplus W) \subseteq T(V\oplus W)$ is generated by 
$I(V), I(W)$ and $x_{i}y_{j} - \chi^{-1}_{j}(K_{i})y_{j} x_{i}$ 
for all $1\leq i,j\leq\theta$, see \cite[Rmk. 1.10]{ARS}. In particular, 
we have that $\widetilde{\U}(\D_{red}, 0) \simeq \toba(V\oplus W)\#\ZZ^{2\theta}$.   
\end{obs}

We end this subsection with a result that translates the notion of twist-equivalence of matrices
of diagonal braidings \cite[Prop. 2.2]{AS3}, \cite[Thm. 2.1]{Ro} to cocycle 
deformations. It states that, under some assumptions, 
these families of pointed Hopf algebras depend only on one parameter 
for each connected component, up to cocycle deformations. 

\begin{theorem}\label{thm:cocycle-def-pre-Nichols}
Assume $\Ga$ is a free abelian group of rank $2\theta$ 
with generators $(L_{i})_{1\leq i\leq \theta}$, $(K_{i})_{1\leq i\leq \theta}$.
Let $\D_{red}$ 
be a reduced YD-datum of Cartan type and $\ell=(\ell_{i}) $ a linking parameter.
If $1\neq q_{ii} $ is a positive real number for all $1\leq i\leq \theta$, then 
$\widetilde{\U}(\D_{red}, \ell)$ is a cocycle deformation of a Hopf algebra 
$\widetilde{\U}(\D_{q}, \ell)$ associated with a reduced YD-datum $\D_{q}$ of DJ-type.
\end{theorem}

\pf By \cite[Thm. 2.1]{Ro} we have that  the Cartan matrix $(a_{ij})_{1\leq i,j\leq \theta}$  
is symmetrizable, with
symmetrizing diagonal matrix $(d_{i})_{1\leq i \leq \theta}$ 
and there is a collection of positive numbers 
$(q_{I})_{I \in \mathcal{X}}$ such
that $(q_{ij})$ is twist-equivalent to $(\hat{q}_{ij})$, where
$\hat{q}_{ij} = q_{I}^{ d_{i} a_{ij}}$
for all $i, j \in I$. 

If we order the group generators 
by $L_{1},\ldots, L_{\theta}, K_{1},\ldots, K_{\theta}$ and take
the corresponding characters $\chi_{1}^{-1},\ldots, \chi_{\theta}^{-1}, \chi_{1},\ldots, \chi_{\theta}$, 
the matrix of the braiding in $V\oplus W$ is given by
$$
p_{ij} = \begin{cases}
             \begin{array}{ll}
              q_{ji}^{-1} &\text{ if }1\leq i,j\leq \theta,\\
              q_{ij}^{-1} &\text{ if }1\leq i \leq \theta,\ \theta+1\leq j\leq 2\theta,\\
                            q_{ji} &\text{ if }\theta+1\leq i\leq 2\theta,\ 1\leq j\leq \theta,\\
                                          q_{ij} &\text{ if }\theta+1\leq i, j\leq 2\theta.
              \end{array}
            \end{cases}$$
Let $\D_{q} = 
\D(\Gamma, (L_{i})_{1\leq i\leq\theta} ,  (K_{i})_{1\leq i\leq\theta},
(q_{I})_{I\in \mathcal{X}},(a_{ij})_{1\leq i,j \leq \theta})$ be the reduced YD-datum 
of DJ-type associated with $(\hat{q}_{ij}) $. Denote by 
$\widehat{V}, \widehat{W}$ the braided vector spaces associated with this datum and by  
$(\hat{p}_{ij})_{i,j }$ the matrix of the braiding in $\widehat{V}\oplus \widehat{W}$. 
Let $\widetilde{\U}(\D_{q}, \ell)$ be the corresponding pointed Hopf 
algebra. 

If we set $g_{i}=L_{i}$ and $g_{i+\theta} = K_{i}$
for all $1\leq i\leq \theta$, then by \cite[Prop. 2.2]{AS3}, 
the map $\sigma: \Ga \times \Ga \to \k^{\times}$ given by
$$ \sigma(g_{i},g_{j}) = \begin{cases}
\begin{array}{ll}
                          \hat{p}_{ij}p_{ij}^{-1} &\text{ if }i\leq j,\\
1 & \text{ otherwise,}
\end{array}                         
\end{cases}
$$
is a group 2-cocycle.
Denote by $\tilde{\sigma}\in \Z^{2}(T(V\oplus W)\#\k  \Ga, \k )$ 
the Hopf 2-cocycle induced by $\sigma$. Then
$\tilde{\sigma}$ induces a Hopf 2-cocycle on $\widetilde{\U}(\D_{red}, \ell)$
and we have that
$\widetilde{\U}(\D_{red}, \ell)_{\tilde{\sigma}} \simeq \widetilde{\U}(\D_{q}, \ell)$. Indeed,
by Proposition \ref{prop:cocycle-lifting-bos} and the proof of \cite[Prop. 3.9]{AS2}
we have that $(T(V\oplus W)\#\k  \Ga)_{\tilde{\sigma}} =T(V\oplus W)_{\sigma}\#\k  \Ga =
T(\hat{V}\oplus \hat{W})\#\k  \Ga$.
For example, for $i\leq j \in I$
$$ K_{i}\cdot_{\sigma} x_{j} = \sigma(K_{i},K_{j})\sigma^{-1}(K_{j},K_{i}) K_{i}\cdot x_{j}
= \hat{q}_{ij}q_{ij}^{-1}q_{ij} x_{j} = \hat{q}_{ij} x_{j} = q_{I}^{d_{i}a_{ij}}x_{j}.$$
Let $J$ be the ideal of $T(V\oplus W)$ generated by the elements 
$x_{i}y_{j} - q_{ij}^{-1}y_{j} x_{i} 
- \delta_{ij} \ell_{i} (K_{i}L_{i} - 1)$ for all $1\leq i,j\leq\theta$. To prove the claim it suffices
to show that the corresponding ideal in $T(V\oplus W)_{\sigma}$ coincides with 
the ideal $\hat{J}$ generated by the elements $x_{i}y_{j} - \hat{q}_{ij}^{-1}y_{j} x_{i} 
- \delta_{ij} \ell_{i} (K_{i}L_{i} - 1)$. But by Proposition \ref{prop:cocycle-lifting-bos} and
the definition of $\sigma$ we have for all $1\leq i,j\leq \theta$ that 
\begin{align*}
x_{i}\cdot_{\sigma}y_{j} - q_{ij}^{-1}y_{j}\cdot_{\sigma} x_{i} 
- \delta_{ij} \ell_{i} (K_{i}\cdot_{\sigma} L_{i} - 1) & = 
\sigma(K_{i},L_{j})x_{i}y_{j} - q_{ij}^{-1}\sigma(L_{j},K_{i}) y_{j}x_{i} 
- \delta_{ij} \ell_{i} (K_{i}L_{i} - 1)\\
& =x_{i}y_{j} - q_{ij}^{-1}(\hat{q}_{ij}q_{ij}^{-1})^{-1}y_{j} x_{i} 
- \delta_{ij} \ell_{i} (K_{i}L_{i} - 1)\\
& =x_{i}y_{j} - \hat{q}_{ij}^{-1}y_{j}x_{i} 
- \delta_{ij} \ell_{i} (K_{i}L_{i} - 1).
\end{align*}
\epf

\subsection{Multiparameter quantum groups as quotients of bosonizations of pre-Nichols algebras.}
\label{subsubsec:explicitex}
In this subsection we show how the
multiparameter quantum groups 
$\QEA$, associated with
a symmetrizable generalized Cartan matrix, introduced
by Hu, Pei and Rosso \cite{HPR},
can be described using reduced data. 
These multiparameter quantum groups contain in a unified way
families of quantum groups introduced by other authors, see \cite{HPR} and
references therein.
Note that in \cite{AS3}, Andruskiewitsch and Schneider 
characterized all pointed Hopf algebras that
can be constructed using a generic datum of finite 
Cartan type for a free group of finite rank. 

Let $\lieg_{A}$ be a symmetrizable Kac-Moody algebra with $A=(a_{ij})_{i,j\in I}$
the associated generalized Cartan matrix, with $I$ a 
finite set. Let $d_{i}$ be relatively prime 
positive integers such that $d_{i}a_{ij} = d_{j}a_{ji}$ for all 
$i,j\in I$. 
Let $\Phi$ be a 
finite root system with 
$\Pi=\{\alpha_{i}:\ i\in I\}$ a set of simple roots,
$Q=\bigoplus_{i\in I}\ZZ \alpha_{i}$ the root lattice,
$\Phi^{+}$ the set of positive roots with respect to $\Pi$ and
$ Q^{+}=\bigoplus_{i\in I}\ZZ_{+} \alpha_{i}$ the positive
root lattice. Let $\textbf{q}=(q_{ij})_{i,j\in I}$ with $q_{ij} \in \k^{\times}$ and 
$q_{ii}\neq 1$ for all $i,j\in I$ satisfying 
\begin{equation}\label{eq:hpr-cartan}
q_{ij}q_{ji} = q_{ii}^{a_{ij}}\qquad \text{for all }i,j\in I.
\end{equation}

\begin{definition}\cite[Def. 7]{HPR}\label{def:multiqgrouphpr}
Let $\QEA$ be the unital associative algebra over $\k $ 
generated by elements $e_{i}, f_{i}, \omega_{i}^{\pm 1}$ and ${\omega'}^{\pm 1}_{i}$ 
with  $i \in I$ satisfying the following relations:
\begin{align*}
(R1)\qquad & \omega_{i}^{\pm 1}{\omega'}^{\pm 1}_{j}  
= {\omega'}^{\pm 1}_{j}{\omega}^{\pm 1}_{i}, 
\qquad \omega_{i}^{\pm 1}{\omega}^{\mp 1}_{i} 
= {\omega'}^{\pm 1}_{i}{\omega'}^{\mp 1}_{i}
= 1,\\
(R2)\qquad & \omega_{i}^{\pm 1}{\omega}^{\pm 1}_{j}  
= {\omega}^{\pm 1}_{j}{\omega}^{\pm 1}_{i}, 
\qquad {\omega'}_{i}^{\pm 1}{\omega'}^{\pm 1}_{j} 
= {\omega'}^{\pm 1}_{j}{\omega'}^{\pm 1}_{i},\\
(R3)\qquad & \omega_{i} e_{j} \omega_{i}^{-1}= q_{ij} e_{j},\qquad
\omega'_{i} e_{j} {\omega'}_{i}^{-1}= q_{ji}^{-1} e_{j},\\
(R4)\qquad &\omega_{i} f_{j} \omega_{i}^{-1}= q_{ij}^{-1} f_{j},\qquad
\omega'_{i} f_{j} {\omega'}_{i}^{-1}= q_{ji}f_{j},\\
(R5)\qquad & [e_{i} , f_{j} ]  = 
\delta_{i,j}\frac{q_{ii}}{q_{ii}-1} (\omega_{i} - \omega'_{i} ),\\
(R6)\qquad & \sum_{k=0}^{1-a_{ij}} (-1)^{k} 
\binom{1-a_{ij}}{k}_{q_{ii}}q_{ii}^{\frac{k(k-1)}{2}}
q_{ij}^{k}e_{i}^{1-a_{ij}-k} e_{j}e_{i}^{k} = 0 \qquad (i\neq j),\\
(R7)\qquad & \sum_{k=0}^{1-a_{ij}} (-1)^{k} 
\binom{1-a_{ij}}{k}_{q_{ii}}q_{ii}^{\frac{k(k-1)}{2}}
q_{ij}^{k} f_{i}^{k} f_{j}f_{i}^{1-a_{ij}-k}= 0 \qquad (i\neq j). 
\end{align*}
\end{definition}

 $\QEA$ is a Hopf algebra with its coproduct, 
counit and antipode 
determined for all $i, j \in I$ by:
\begin{align*}
\com(e_{i}) & = e_{i}\ot 1 + \omega_{i}\ot e_{i}, &
\eps(e_{i}) & = 0,& \SS(e_{i}) = -\omega_{i}^{-1} e_{i},\\
\com(f_{i}) & = f_{i}\ot \omega'_{i} + 1\ot f_{i}, &
\eps(f_{i}) & = 0,& \SS(f_{i}) = - f_{i}{\omega'_{i}}^{-1},\\
\com(\omega_{i}^{\pm 1}) & = \omega_{i}^{\pm 1}\ot \omega_{i}^{\pm 1},&
\eps(\omega_{i}^{\pm 1}) & = 1, 
& \SS(\omega_{i}^{\pm 1}) = \omega_{i}^{\mp 1}, \\
\com({\omega'}_{i}^{\pm 1}) & = {\omega'}_{i}^{\pm 1}\ot {\omega'}_{i}^{\pm 1}, &
\eps({\omega'}_{i}^{\pm 1})& = 1, 
& \SS({\omega'}_{i}^{\pm 1}) = {\omega'}_{i}^{\mp 1}. 
\end{align*}

Next we prove that this quantum group can be described using reduced data.

\subsection*{Definition of $\widetilde{\U}(\D_{red}, \ell)$}
Let 
\begin{itemize}
 \item $\theta = |I|$,
 \item $\Ga = \ZZ^{2|I|}$ and denote $K_{i}, L_{i}$ 
with $i\in I$ two (commuting) generators,
 
 \item $\chi_{i}\in \widehat{\Ga}$ given by $\chi_{i}(K_{j}) = q_{ji}$ and 
$\chi_{i}(L_{j}) = q_{ij}$ for all $i,j\in I$.
\end{itemize}
In particular, we have that $\chi_{i}(L_{j}) = \chi_{j}(K_{i})$ for all $i,j\in I$.
Since by assumption 
$K_{i}L_{i} \neq  1$ and by \eqref{eq:hpr-cartan},
$q_{ij}q_{ji} = q_{ii}^{a_{ij}} $ with $q_{ii}\neq 1$,
we have 
that $\D_{red} = \D(\Ga, (K_{i}), (L_{i}), (\chi_{i}), (a_{ij}) )$ 
is a reduced YD-datum of 
Cartan type. 

Let $V$,
$W$ be the vector spaces linearly generated by the elements
$x_{i}$ and $y_{i}$
for all $1\leq i\leq \theta$. Following the definition of reduced data, 
both have a Yetter-Drinfeld module 
structure. In this case, it is given for all $i,j\in I$ by 
\begin{align*}
\delta(x_{j}) & = K_{j} \ot x_{j}, & K_{i}\cdot x_{j}  &= \chi_{j}(K_{i})x_{j}=
q_{ij} x_{j},
& L_{i}\cdot x_{j} &=
q_{ji}x_{j},
\\
\delta(y_{j}) & = L_{j} \ot y_{j}, & 
K_{i}\cdot y_{j} &= \chi_{j}^{-1}(K_{i})y_{j}=
q_{ij}^{-1}y_{j},
 & L_{i}\cdot y_{j} &=
q_{ji}^{-1}y_{j}.
\end{align*}
 
Recall that for $\ell  = (\ell_{i})_{1\leq i \leq \theta}$ with $\ell_{i} \in \k^{\times}$,
the pointed Hopf algebra $\widetilde{\U}(\D_{red}, \ell)$ associated
with these data is given by 
the quotient Hopf algebra of the bosonization $R (\D) \# \k  \ZZ^{2\theta}$ 
modulo the
ideal generated by
\begin{align*}
x_{i}y_{j} - q_{ij}^{-1}y_{j} x_{i} 
- \delta_{ij} \ell_{i} (K_{i}L_{i} - 1)\qquad\text{for all } i,j\in I. 
\end{align*}

In particular, in $\widetilde{\U} (\D_{red},\ell)$,
$x_{i}$ is a $(1,K_{i})$-primitive and $y_{i} $
is a $(1,L_{i})$-primitive. Indeed, for $x_{i} \in R (\D) \# \k  \ZZ^{2\theta}$ we have
$\com_{R}(x_{i}) = x_{i}^{(1)}\ot x_{i}^{(2)} = x_{i}\ot 1 + 1 \ot x_{i}$ and 
$$
\com(x_{i})  = x_{i}^{(1)}\#(x_{i}^{(2)})_{(-1)}\ot (x_{i}^{(2)})_{(0)}\#1
= (x_{i} \# 1)\ot (1\#1) + (1\#K_{i})\ot (x_{i}\#1)
= x_{i}\ot 1 +  K_{i}\ot x_{i}. 
$$
 
\begin{theorem}\label{thm:qersg=ured}
$\QEA \simeq \widetilde{\U} (\D_{red},\ell)$ with 
$\ell_{i} = \frac{q_{ii}}{q_{ii}-1}$ for all 
$1\leq i\leq \theta$. 
\end{theorem}

\pf 
Let $\varphi: \QEA \to \widetilde{\U}(\D_{red},\ell)$ be the algebra map 
defined by 
$$\varphi(\omega_{i}) = K_{i},\quad 
\varphi(\omega'_{i}) = L_{i}^{-1},\quad
\varphi(e_{i}) = x_{i},\quad
\varphi(f_{i}) = y_{i}L_{i}^{-1} \qquad \text{ for all }1\leq i\leq \theta.$$ 
The map $\varphi$ is an epimorphism, if it is well-defined.
To prove that it is well-defined, we show that the relations in $\QEA$
are mapped to $0$ by $\varphi$. First, notice that the action of
$K_{i}$ and $L_{i}$ on $x_{j}$ and $y_{j}$ yields a commutation
relation in $T(V\oplus W)\# \k  \ZZ^{2\theta}$; for example,  
$(1\# K_{i})(x_{j}\# 1) (1\#K_{i}^{-1}) = 
[(K_{i}\cdot x_{j} \# K_{i} ) (1\#K_{i}^{-1})] = 
K_{i}\cdot x_{j} \# 1
$.
Clearly, we need to verify only relations $(R3)-(R7)$.
For $(R3)$, we have 
\begin{align*}
\varphi(\omega_{i}e_{j}\omega_{i}^{-1} - 
q_{ij}e_{j}) & = 
K_{i}x_{j}K_{i}^{-1} - q_{ij}x_{j} = 
K_{i}\cdot x_{j} - q_{ij}x_{j} = 0,\\
\varphi(\omega'_{i}e_{j}{\omega'}_{i}^{-1} - 
q_{ji}^{-1}e_{j}) & = 
L_{i}^{-1}x_{j}L_{i} - 
q_{ji}^{-1}x_{j} = 
L_{i}^{-1}\cdot x_{j} - q_{ji}^{-1}x_{j}
= 0. 
\end{align*}
The proof for $(R4)$ follows the same lines. Since 
$\D$ is a reduced datum, for $(R5)$ we have
\begin{align*}
\varphi([e_{i},f_{j}]) & = x_{i}y_{j}L_{j}^{-1} - y_{j}L_{j}^{-1}x_{i} 
=  x_{i}y_{j}L_{j}^{-1} - \chi_{i}(L_{j}^{-1})y_{j}x_{i}L_{j}^{-1}
=  x_{i}y_{j}L_{j}^{-1} - q_{ij}^{-1}y_{j}x_{i}L_{j}^{-1}\\
&=  (x_{i}y_{j} - q_{ij}^{-1}y_{j}x_{i})L_{j}^{-1}
= \delta_{ij}\ell_{i}(K_{i}L_{i} -1)L_{j}^{-1}
 = \delta_{ij}\frac{q_{ii}}{q_{ii}-1}(K_{i} - L_{i}^{-1}) \\
& =  \varphi\left(\delta_{ij}\frac{q_{ii}}{q_{ii}-1}(\omega_{i}- \omega'_{i})\right).
\end{align*}
 
To verify $(R6)$ and $(R7)$ one notes that their images under $\varphi$
are the quantum Serre relations in $T(V\oplus W)$, e.~g.~ 
$\ad_{c}(x_{i})^{1-a_{ij}}(x_{j})$ is the image
of the left hand side of $(R6)$. Indeed, since
\begin{align*}
\ad_{c}(x_{i})(x_{j}) & = x_{i}x_{j} - [(x_{i})_{-1}\cdot x_{j}] (x_{i})_{0}
= x_{i}x_{j} - (K_{i}\cdot x_{j}) x_{i}
 = x_{i}x_{j} - q_{ij}x_{j} x_{i}.
 \end{align*}
one proves by induction that 
\begin{equation}\label{eq:qserrexi}
\ad_{c}(x_{i})^{n}(x_{j})=  
 \sum_{k=0}^{n} (-1)^{k} 
\binom{n}{k}_{q_{ii}}q_{ii}^{\frac{k(k-1)}{2}}
q_{ij}^{k}x_{i}^{n-k} x_{j}x_{i}^{k} \quad \text{ for all }n\in \NN.
 \end{equation}
Assuming that the equality holds for $n\in \NN$ 
and using \eqref{eq:qbinom} we have
\begin{align*}
& \ad_{c}(x_{i})^{n+1}(x_{j}) = \ad_{c}(x_{i})(\ad_{c}(x_{i})^{n}(x_{j})) 
= \ad_{c}(x_{i})\left(\sum_{k=0}^{n} (-1)^{k} 
\binom{n}{k}_{q_{ii}}q_{ii}^{\frac{k(k-1)}{2}}
q_{ij}^{k}x_{i}^{n-k} x_{j}x_{i}^{k} \right) 
\\
&\quad = \sum_{k=0}^{n} (-1)^{k} 
\binom{n}{k}_{q_{ii}}q_{ii}^{\frac{k(k-1)}{2}}
q_{ij}^{k}[x_{i}^{n+1-k} x_{j}x_{i}^{k} - K_{i}
\cdot(x_{i}^{n-k} x_{j}x_{i}^{k})x_{i}]
\\
&\quad = \sum_{k=0}^{n} (-1)^{k} 
\binom{n}{k}_{q_{ii}}q_{ii}^{\frac{k(k-1)}{2}}
q_{ij}^{k}[x_{i}^{n+1-k} x_{j}x_{i}^{k} - q_{ii}^{n}q_{ij}
x_{i}^{n-k} x_{j}x_{i}^{k+1}]
\\
&\quad = \sum_{k=0}^{n} (-1)^{k} 
\binom{n}{k}_{q_{ii}}q_{ii}^{\frac{k(k-1)}{2}}
q_{ij}^{k}x_{i}^{n+1-k} x_{j}x_{i}^{k} + 
\sum_{k=0}^{n} (-1)^{k+1} 
\binom{n}{k}_{q_{ii}}q_{ii}^{\frac{k(k-1)}{2}}
q_{ij}^{k+1}q_{ii}^{n}
x_{i}^{n-k} x_{j}x_{i}^{k+1}
\\
&\quad = x_{i}^{n+1}x_{j} + (-1)^{n+1} 
q_{ii}^{\frac{n(n+1)}{2}}q_{ij}^{n+1}x_{j}x_{i}^{n+1} +\\
&\quad\quad +  \sum_{k=1}^{n} (-1)^{k} q_{ij}^{k}\left[
\binom{n}{k}_{q_{ii}}q_{ii}^{\frac{k(k-1)}{2}} + 
\binom{n}{k-1}_{q_{ii}}q_{ii}^{\frac{(k-1)(k-2)}{2}}q_{ii}^{n}\right]
x_{i}^{n+1-k} x_{j}x_{i}^{k} 
\\
&\quad = x_{i}^{n+1}x_{j} + (-1)^{n+1} 
q_{ii}^{\frac{n(n+1)}{2}}q_{ij}^{n+1}x_{j}x_{i}^{n+1} +\\
&\quad\quad +  \sum_{k=1}^{n} (-1)^{k} q_{ij}^{k}q_{ii}^{\frac{k(k-1)}{2}}
\left[\binom{n}{k}_{q_{ii}} + 
\binom{n}{k-1}_{q_{ii}}q_{ii}^{n+1-k}\right]
x_{i}^{n+1-k} x_{j}x_{i}^{k}
\end{align*}
\begin{align*}
&\quad = x_{i}^{n+1}x_{j} + (-1)^{n+1} 
q_{ii}^{\frac{n(n+1)}{2}}q_{ij}^{n+1}x_{j}x_{i}^{n+1} + 
\sum_{k=1}^{n} (-1)^{k} q_{ij}^{k}q_{ii}^{\frac{k(k-1)}{2}}
\binom{n+1}{k}_{q_{ii}}
x_{i}^{n+1-k} x_{j}x_{i}^{k}\\
&\quad =\sum_{k=0}^{n+1} (-1)^{k} q_{ij}^{k}q_{ii}^{\frac{k(k-1)}{2}}
\binom{n+1}{k}_{q_{ii}}
x_{i}^{n+1-k} x_{j}x_{i}^{k}.
\end{align*}
Since $\ad_{c}(x_{i})^{1-a_{ij}}(x_{j}) = 0$ in $R(\D)$, and
$$\varphi\left( 
\sum_{k=0}^{1-a_{ij}} (-1)^{k} 
\binom{1-a_{ij}}{k}_{q_{ii}}q_{ii}^{\frac{k(k-1)}{2}}
q_{ij}^{k}e_{i}^{1-a_{ij}-k} e_{j}e_{i}^{k}\right)=\ad_{c}(x_{i})^{1-a_{ij}}(x_{j}),$$
the assertion about $(R6)$ follows.
Analogously, 
\begin{align*}
\ad_{c}(y_{i})(y_{j}) & =
y_{i}y_{j} - (L_{i}\cdot y_{j}) y_{i}
= y_{i}y_{j} - q_{ji}^{-1}y_{j} y_{i} = -q_{ji}^{-1}(y_{j}y_{i} - q_{ji}y_{i}y_{j}).
 \end{align*}
and one may prove by induction that 
\begin{equation}\label{eq:qserreyi}
\ad_{c}(y_{i})^{n}(y_{j}) =  
 (-1)^{n}q_{ji}^{-n}q_{ii}^{-\frac{n(n-1)}{2}}\left(\sum_{k=0}^{n} (-1)^{k} 
\binom{n}{k}_{q_{ii}}q_{ii}^{\frac{k(k-1)}{2}}
q_{ji}^{k}y_{i}^{k} y_{j}y_{i}^{n-k}\right).
 \end{equation}
Hence,
\begin{align*}
& \varphi\left(\sum_{k=0}^{1-a_{ij}} (-1)^{k} 
\binom{1-a_{ij}}{k}_{q_{ii}}q_{ii}^{\frac{k(k-1)}{2}}
q_{ij}^{k} f_{i}^{k} f_{j}f_{i}^{1-a_{ij}-k}\right) = \\
& \quad = \sum_{k=0}^{1-a_{ij}} (-1)^{k} 
\binom{1-a_{ij}}{k}_{q_{ii}}q_{ii}^{\frac{k(k-1)}{2}}
q_{ij}^{k} (y_{i}L_{i}^{-1})^{k} (y_{j}
L_{j}^{-1})(y_{i}L_{i}^{-1})^{1-a_{ij}-k}
\\
& \quad = \sum_{k=0}^{1-a_{ij}} (-1)^{k} 
\binom{1-a_{ij}}{k}_{q_{ii}}q_{ii}^{\frac{k(k-1)}{2}}
q_{ij}^{k} q_{ii}^{\frac{k(k-1)}{2}}
q_{ii}^{\frac{(1-a_{ij}-k)(-a_{ij}-k)}{2}}
y_{i}^{k}L_{i}^{-k} y_{j}L_{j}^{-1}
y_{i}^{1-a_{ij}-k}L_{i}^{-1+a_{ij}+k}
\end{align*}
\begin{align*}
& \quad = (\sum_{k=0}^{1-a_{ij}} (-1)^{k} 
\binom{1-a_{ij}}{k}_{q_{ii}}q_{ii}^{k(k-1) }
q_{ii}^{\frac{(1-a_{ij}-k)(-a_{ij}-k)}{2}}
q_{ji}^{k}q_{ii}^{k(1-a_{ij}-k)}q_{ij}^{1-a_{ij}}
y_{i}^{k}y_{j}y_{i}^{1-a_{ij}-k})L_{j}^{-1}
L_{i}^{-1+a_{ij}}
\end{align*}
But
\begin{align*}
& {\frac{k(k-1)}{2}+
\frac{(1-a_{ij}-k)(-a_{ij}-k)}{2} + k(1-a_{ij}-k)}
= \\
& \qquad = {\frac{k(k-1)}{2}+
\frac{-a_{ij}(1-a_{ij}-k)-k(1-a_{ij}-k))}{2} + k(1-a_{ij}-k)}\\
& \qquad = {\frac{k(k-1)}{2}+
\frac{-a_{ij}(1-a_{ij}-k)+k(1-a_{ij}-k)}{2}}\\
& \qquad= \frac{1}{2}[k^{2} - k - a_{ij} + a_{ij}^{2} + a_{ij}k
+ k -ka_{ij} - k^{2}] = \frac{1}{2} a_{ij}(a_{ij}-1). 
\end{align*}
Thus, $\varphi$ of $(R7)$ equals 
\begin{align*}
\left(\sum_{k=0}^{1-a_{ij}} (-1)^{k} 
\binom{1-a_{ij}}{k}_{q_{ii}}q_{ii}^{\frac{k(k-1)}{2}} q_{ji}^{k}
q_{ii}^{\frac{a_{ij}(a_{ij}-1)}{2}}q_{ij}^{1-a_{ij}}
y_{i}^{k}y_{j}y_{i}^{1-a_{ij}-k}\right)L_{j}^{-1}
L_{i}^{-1+a_{ij}}\\
 \qquad = 
q_{ij}^{1-a_{ij}}
(-1)^{1-a_{ij}}q_{ji}^{1-a_{ij}}
\ad_{c}(y_{i})^{1-a_{ij}}(y_{j})
L_{j}^{-1}
L_{i}^{-1+a_{ij}}
\end{align*}
Since $\ad_{c}(y_{i})^{1-a_{ij}}(y_{j}) = 0$ in $R(\D)$, the claim about
$(R7)$ follows.

Hence, $\varphi$ is a well-defined algebra map. 
Moreover, it is a Hopf algebra map, since 
$\omega_{i}$, $\omega'_{i}$ and
$K_{i}$, $L_{i}$
are grouplike elements, 
$e_{i}$ is $(1,\omega_{i})$-primitive and 
$f_{i}$ is $(\omega'_{i},1)$-primitive, and
the elements $x_{i}$ and $y_{i}$ are
$(1,K_{i})$-primitive and $(1, L_{i})$-primitive, respectively,
for all $i \in I$.

Now we show $\varphi$ is an isomorphism.
Let $\widetilde{\psi}: R (\D) \# \k  \ZZ^{2\theta} \to \QEA  $ 
be the algebra map 
given by 
$$\widetilde{\psi}(1\#K_{i}) = \omega_{i},\quad 
\widetilde{\psi}(1\#L_{i}) = {\omega'}^{-1}_{i},\quad
\widetilde{\psi}(x_{i}\#1) = e_{i},\quad
\widetilde{\psi}(y_{i}\#1) = f_{i}{\omega'}_{i}^{-1},$$ 
for all $i \in \I$. 
Again, $\widetilde{\psi}$ is clearly a 
Hopf algebra epimorphism, if it is well-defined. 
For example, it preserves the algebra structure, if $i \in \I$ we have
\begin{align*}
\widetilde{\psi}((1\#K_{i})(x_{j}\# 1)(1\#K_{i}^{-1})) = 
\omega_{i}e_{j}\omega_{i}^{-1} = q_{ij}e_{j}=
\chi_{j}(K_{i}) e_{j} = \widetilde{\psi}(K_{i}\cdot x_{j}\# 1),  
\end{align*}
and  the coalgebra structure with 
$\varepsilon(\omega_{i}) = 1 = \varepsilon(1\#K_{i})$,
$\varepsilon(\omega'_{i}) = 1 = \varepsilon(1\#L_{i})$,
$\varepsilon(e_{i}) = 0 = \varepsilon(x_{i}\#1)$,
$\varepsilon(f_{i}{\omega'}_{i}^{-1}) = 1 = \varepsilon(y_{j}\#1)$ and
\begin{align*}
\com(\widetilde{\psi}(y_{i}\#1) )  & =
\com( f_{i}{\omega'}_{i}^{-1})  = 
f_{i}{\omega'}_{i}^{-1} \otimes 1 +
 {\omega'}_{i}^{-1}\otimes f_{i}{\omega'}_{i}^{-1}\\ 
& =  \widetilde{\psi}(y_{i}\# 1)\otimes 1 + \widetilde{\psi}(1\#L_{i}) \otimes 
\widetilde{\psi}(y_{i}\#1) = (\widetilde{\psi}\otimes \widetilde{\psi})\com(y_{i}\#1).
\end{align*}
To see that $\widetilde{\psi}$ is indeed well-defined, we have to check
that the quantum Serre relations are mapped to $0$.
For $i\neq j \in I$ we have by \eqref{eq:qserrexi} and $(R6)$ that 
$$
 \widetilde{\psi}(\ad_{c}(x_{i})^{1-a_{ij}}(x_{j})) = 
 \sum_{k=0}^{1-a_{ij}} (-1)^{k} 
\binom{1-a_{ij}}{k}_{q_{ii}}q_{ii}^{\frac{k(k-1)}{2}}
q_{ij}^{k}e_{i}^{1-a_{ij}-k} e_{j}e_{i}^{k} =0.
$$
Analogously, by \eqref{eq:qserreyi}, $(R7)$ and the same 
calculation with the exponents as above we have
\begin{align*}
& \widetilde{\psi}(\ad_{c}(y_{i})^{1-a_{ij}}(y_{j})) = 
(-1)^{1-a_{ij}}q_{ji}^{-1+a_{ij}}q_{ii}^{\frac{a_{ij}(1-a_{ij})}{2}}\cdot \\
& \quad\qquad \cdot \left(\sum_{k=0}^{1-a_{ij}} (-1)^{k} 
\binom{1-a_{ij}}{k}_{q_{ii}}q_{ii}^{\frac{k(k-1)}{2}}
q_{ji}^{k}(f_{i}{\omega'}^{-1}_{i})^{k} 
(f_{j}{\omega'}^{-1}_{j})(f_{i}{\omega'}^{-1}_{i})^{1-a_{ij}-k}\right)\\
& \quad = (-1)^{1-a_{ij}}q_{ji}^{-1+a_{ij}} q_{ii}^{\frac{a_{ij}(1-a_{ij})}{2}}
\sum_{k=0}^{1-a_{ij}} (-1)^{k} 
\binom{1-a_{ij}}{k}_{q_{ii}}q_{ii}^{\frac{k(k-1)}{2}}
q_{ji}^{k}\cdot \\
& \quad\qquad
\cdot q_{ii}^{\frac{k(k-1)}{2}+\frac{(1-a_{ij}-k)(1-a_{ij})}{2}}
f_{i}^{k}{\omega'}^{-k}_{i}  f_{j}{\omega'}^{-1}_{j}f_{i}^{1-a_{ij}-k}
{\omega'}^{-1+a_{ij}+k}_{i})\\
& \quad = (-1)^{1-a_{ij}}q_{ji}^{-1+a_{ij}} q_{ii}^{\frac{a_{ij}(1-a_{ij})}{2}}
\sum_{k=0}^{1-a_{ij}} (-1)^{k} 
\binom{1-a_{ij}}{k}_{q_{ii}}q_{ii}^{\frac{k(k-1)}{2}}
q_{ji}^{k}\cdot \\
& \quad\qquad 
\cdot q_{ii}^{-\frac{k(k-1)}{2}-\frac{(1-a_{ij}-k)(-a_{ij}-k)}{2}}
q_{ji}^{-k}q_{ij}^{-1+a_{ij}+k}q_{ii}^{-k(1-a_{ij}-k)}
f_{i}^{k}  f_{j}f_{i}^{1-a_{ij}-k}
{\omega'}_{j}{\omega'}^{-1+a_{ij}}_{i}\\
& \quad = (-1)^{1-a_{ij}}q_{ji}^{-1+a_{ij}} q_{ii}^{\frac{a_{ij}(1-a_{ij})}{2}}
\sum_{k=0}^{1-a_{ij}} (-1)^{k} 
\binom{1-a_{ij}}{k}_{q_{ii}}q_{ii}^{\frac{k(k-1)}{2}}q_{ij}^{k}
\cdot \\
& \quad\qquad
\cdot q_{ii}^{-\frac{k(k-1)}{2}-\frac{(1-a_{ij}-k)(-a_{ij}-k)}{2}-k(1-a_{ij}-k) }
q_{ij}^{-1+a_{ij}}
f_{i}^{k}  f_{j}f_{i}^{1-a_{ij}-k}
{\omega'}_{j}{\omega'}^{-1+a_{ij}}_{i}\\
& \quad = (-1)^{1-a_{ij}}q_{ji}^{-1+a_{ij}}
q_{ij}^{-1+a_{ij}}
\left(\sum_{k=0}^{1-a_{ij}} (-1)^{k} 
\binom{1-a_{ij}}{k}_{q_{ii}}q_{ii}^{\frac{k(k-1)}{2}}q_{ij}^{k}
f_{i}^{k}  f_{j}f_{i}^{1-a_{ij}-k}\right)
{\omega'}_{j}{\omega'}^{-1+a_{ij}}_{i}=0.
\end{align*}

Moreover, identifying $x_{i} = x_{i}\# 1$, $y_{i}=y_{i}\#1$ and 
$K_{i}= 1\# K_{i}$, $L_{i}= 1\# L_{i}$, 
by $(R5)$ we have that
\begin{align*}
&\widetilde{\psi}(x_{i}y_{j} - q_{ij}^{-1}y_{j} x_{i}) =
e_{i}f_{j}{\omega'}_{j}^{-1} - 
q_{ij}^{-1}f_{j}{\omega'}_{j}^{-1} e_{i}
= (e_{i}f_{j}- q_{ij}^{-1}q_{ij}
f_{j}e_{i}){\omega'}_{j}^{-1} = 
 [e_{i},f_{j}] {\omega'}_{j}^{-1}=\\
& \qquad =
 \delta_{i,j}\frac{q_{ii}}{q_{ii}-1} (\omega_{i} - {\omega'}_{i}){\omega'}_{j}^{-1}
= \delta_{i,j}\frac{q_{ii}}{q_{ii}-1}  (\omega_{i}{\omega'}_{i}^{-1} - 1)
= \widetilde{\psi}(\delta_{i,j} \ell_{i} (K_{i}L_{i} - 1)),
\end{align*}
for all $i,j \in \I$. Thus $\widetilde{\psi}$ induces a Hopf algebra
epimorphism $\psi:  \widetilde{\U} (\D_{red},\ell) \to \QEA$ such 
$\varphi \circ \psi = \id = \psi \circ \varphi$, implying that $\varphi$
is an isomorphism.  
\epf

\begin{cor}\label{cor:multi-one-cocycle}
Assume $\mathfrak{g}_{A}$ is simple and let $U_{q}(\mathfrak{g}_{A})$ be the 
one-parameter quantum group of Drinfeld-Jimbo type. 
Let $\widetilde{\U} (\D_{q},\ell)$ be the pointed Hopf algebra associated
with the reduced YD-datum of DJ-type 
$\D_{q}= (\Gamma, (L_{i})_{1\leq i\leq\theta} ,  (K_{i})_{1\leq i\leq\theta},
q,(a_{ij})_{1\leq i,j \leq \theta})$ and $\ell$ with $\ell_{i} = \frac{q}{q-1}$
for all $i\in I$, and denote 
$\widehat{\U} (\D_{q},\ell) =  \widetilde{\U} (\D_{q},\ell)/ (K_{i} - L_{i}^{-1})$. Then 
$U_{q}(\mathfrak{g}_{A})\simeq \widehat{\U} (\D_{q},\ell)$.
\end{cor}

\pf
By \cite[Rmk. 9]{HPR}, we know that $U_{q}(\mathfrak{g}_{A}) \simeq 
\QEA / (\omega_{i}'-\omega_{i}^{-1})$ for $q_{ij} = q^{d_{i}a_{ij}}$ for all $i,j\in I$. Moreover,
by Theorem \ref{thm:qersg=ured} we know that $\QEA \simeq \widetilde{\U} (\D_{q},\ell)$
for $\D_{q}$ the reduced YD-datum of DJ-type described above, since
$\mathfrak{g}_{A}$ is connected and the assumption on the $q_{ij}$'s.
The result follows since the isomorphism factors through the quotients. 
\epf

\begin{obs}\label{rmk:one-par-redDJtype}
 With the assumptions above, we have that $\QEA = U_{q,q^{-1}}(\mathfrak{g}_{A})$, the 
special case of the two-parameter quantum group, see 
  \cite[Rmk. 9]{HPR}. Then $U_{q,q^{-1}}(\mathfrak{g}_{A}) \simeq \widetilde{\U} (\D_{q},\ell)$.
\end{obs}

\subsection{Multiparameter quantum groups as cocycle deformations}
In this subsection we apply Theorems \ref{thm:cocycle-def-pre-Nichols} and 
\ref{thm:qersg=ured} to multiparameter quantum groups. 
In particular, we obtain another description of 
\cite[Thm. 28]{HPR}. 
\textit{From now on we assume that $1\neq q_{ii}$ is a positive real 
number for all $i\in I$ and all $I \in \mathcal{X}$.  }

 The existence of all positive roots for every  $ q_{ii} $  
ensures that for each $I\in \mathcal{X}$ and for each $i \in I $  there exists  $q_i \in \k^\times$  
 such that  $ q_{ii} = q_i^{\,2\,d_i}$. Moreover,  we may choose 
  in each connected component $J\in \mathcal{X}$ a unique constant value  $q_i =: q_{J} $ for
  all $i\in J$.
 
Let $\widetilde{\U} (\D_{q},\ell)$ be the pointed Hopf algebra associated
to the reduced YD-datum of DJ-type given by
$\D_{q}= \D(\Gamma, (L_{i})_{1\leq i\leq\theta} ,  (K_{i})_{1\leq i\leq\theta},
(q_{J})_{J\in \mathcal{X}},(a_{ij})_{1\leq i,j \leq \theta})$.

 \begin{cor}\label{cor:multi-two-cocycle}
 There exists a group 2-cocycle $\sigma \in \Z^{2}(\Ga,\k)$
 such that 
$\widetilde{\U} (\D_{q},\ell) \simeq \QEA_{\tilde{\sigma}}$.
In particular, if $\mathfrak{g}_{A}$ is simple and $q_{ii} = q^{2d_{i}}$
for all $i,j\in I$, we have that 
$U_{q,q^{-1}}(\mathfrak{g}_{A}) \simeq \QEA_{\tilde{\sigma}}$.
 \end{cor}

\pf By Theorem \ref{thm:qersg=ured} we know that 
$\QEA \simeq \widetilde{\U} (\D_{red},\ell)$ with 
$\ell_{i} = \frac{q_{ii}}{q_{ii}-1}$ for all 
$i\in I$. Since the braiding is positive and generic, Theorem \ref{thm:cocycle-def-pre-Nichols}
implies that $\QEA$ is a cocycle deformation of a pointed Hopf algebra associated to 
a reduced YD-datum of DJ-type. For, 
the proof of Theorem \ref{thm:cocycle-def-pre-Nichols} gives a group 2-cocycle 
${\sigma}$ such that $\widetilde{\U} (\D_{q},\ell)\simeq  \widetilde{\U} (\D_{red},\ell)_{\tilde{\sigma}}$. 
Taking the corresponding 2-cocycle induced by the isomorphism we have the assertion.

If $\mathfrak{g}_{A}$ is simple and $q_{ij} = q^{d_{i}a_{ij}}$
for all $i,j\in I$, by Remark \ref{rmk:one-par-redDJtype} we have that $\widetilde{\U} (\D_{q},\ell)
\simeq U_{q,q^{-1}}(\mathfrak{g}_{A})$.   
\epf

\begin{obs}
Assume $\mathfrak{g}_{A}$ is simple.
The result above was previously obtained by \cite{HPR} where a Hopf 2-cocycle $\sigma$
is defined in $U_{q,q^{-1}}(\mathfrak{g}_{A})$. We show that this cocycle
comes from a group 2-cocycle on $\Ga$.

First we fix the notation  
   $\omega_\lambda := \prod_{i \in I} \omega_i^{\lambda_i}$  and  
   $\omega'_\lambda := 
   \prod_{i \in I} {\omega'_i}^{\lambda_i}$  for every  $\lambda = 
   \sum_{i \in I} \lambda_i \alpha_i \in Q$.  Similarly, we shall also write
   $$ q_{\mu\nu}:=  {\textstyle \prod\limits_{i,j \in I}} \, q_{ij}^{\, \mu_i \nu_j}   
   \qquad \forall \ \mu = {\textstyle \sum\limits_{i \in I}} \, \mu_i \, \alpha_i \; , \,\; \nu = 
   {\textstyle \sum\limits_{j \in I}} \, \nu_j \, \alpha_j \, \in \, Q$$

 Let $\sigma :  U_{q,q^{-1}}(\mathfrak{g}_{A})\otimes 
 U_{q,q^{-1}}(\mathfrak{g}_{A}) \to \k $  
 be the unique  $\k$-linear form on  $ \QEA $  such that
  $$  \sigma(x,y)  \, := \,  \begin{cases}
           \;\  q^{{{\,1\,} \over {\,2\,}}}_{\mu\,\nu}  &  \;\;  \text{if}  \quad  
           x = \omega_\mu  \text{\;\ or \ }  x = \omega'_\mu \; ,  \quad  y = \omega_\nu  \text{\;\ or \ }  y = \omega'_\nu  \\
           \;\  0  &  \;\;  \text{otherwise}.
                             \end{cases}  $$
Then by \cite[27,28]{HPR}, $\sigma \in \Z^{2}(U_{q,q^{-1}}(\mathfrak{g}_{A}),\k)$ 
and $\QEA\simeq U_{q,q^{-1}}(\mathfrak{g}_{A})_{\sigma} $.

On the other hand, we know that $ U_{q,q^{-1}}(\mathfrak{g}_{A})$ is a quotient of a 
bosonization $T(V\oplus W)\# \k \Ga$ with $\Ga = \ZZ^{2\theta}$. As in Remark \ref{rmk:def-ast}, 
we have a $\Ga\times \Ga $
grading on $T(V\oplus W)$
induced by the coaction on the Yetter-Drinfeld module; for example, $\omega_{i}$ have degree
$(\omega_{i},\omega_{i})$, $e_{i}$ has degree $(\omega_{i},1)$
and $f_{i}$ has degree $(1,\omega_{i}'^{-1})$ for all $i\in I$. 
Consider now the group 2-cocycle  $\varphi \in \Z^{2}(\Ga,\k) $  given by
  $\varphi = \sigma|_{\Ga\times \Ga}$, that is,
  $$  \varphi(h,k):= q^{{{\,1\,} \over {\,2\,}}}_{\mu\,\nu} \qquad\text{ if }  \quad  
           h = \omega_\mu  \text{ or }  h = \omega'_\mu \; ,  \quad  k = \omega_\nu  \text{ or }  k = \omega'_\nu,$$
and let $\tilde{\varphi}$ be the 2-cocycle defined on $T(V\oplus W)\# \k \Ga$.
Since the group is abelian and 
$e_{i}\cdot_{\tilde{\varphi}}f_{j} = e_{i}f_{j}$ for all $i,j\in I$, 
we have that $e_{i}\cdot_{\tilde{\varphi}}f_{j} - f_{j}\cdot_{\tilde{\varphi}}e_{i}
= [e_{i}, f_{j}]$ and consequently $ \tilde{\varphi}$ defines 
a Hopf 2-cocycle on $U_{q,q^{-1}}(\mathfrak{g}_{A})$.
Since 
$\sigma(x,y) = 0 = \eps(x)\eps(y)$ if $x,y\notin \Ga$, it follows that 
$\sigma = \tilde{\varphi}$
and whence 
$\QEA =U_{q,q^{-1}}(\mathfrak{g}_{A})_{\sigma} $.
\end{obs}

\section*{Acknowledgements}
Research of this paper was begun when the author was visiting 
Universit\`a di Roma ``Tor Vergata" under 
the support of the GNSAGA and CONICET. He thanks 
F. Gavarini and the people of the Math. Department for the warm
hospitality.  These notes are intended to
contribute to a joint project with F. Gavarini on multiparameter 
quantum groups; the author wishes to thank all conversations and
comments which helped to improve the paper. 
He also wishes to thank A. Garc\'ia Iglesias and
L. Vendramin for interesting discussions.

\bibliographystyle{amsbeta}

\end{document}